\documentclass[final]{elsarticle}
\usepackage[textsize=small]{todonotes}
\reversemarginpar

\journal{Journal of Computational Physics}

\usepackage{amsmath}
\usepackage{amssymb}
\usepackage{afterpage,tabularx}

\usepackage{arydshln}
\usepackage{txfonts}
\usepackage{stmaryrd}
\usepackage{graphicx}
\usepackage{algorithmic}
\usepackage{algorithm}
\usepackage{color}
\usepackage{diagbox,hyperref}
\usepackage{multirow}
\usepackage{rotating}
\usepackage{verbatim}
\usepackage{color}

 \font\msbm=msbm10

\def\dsR{\hbox{{\msbm \char "52}}}

\def\div{\nabla \cdot}

\def\ssl{\mathfrak{B}}
\def\ssr{\mbox{\em l}}

\def\pP{\mathcal{J}}
\def\pL{\mathcal{L}}
\def\pB{\mathfrak{b}}
\def\pC{\mathfrak{c}}

\def\pd{\mathfrak{T}}

\def\pE{\mathbf{E}}

\def\ppE{\tilde{\mathbf{E}}}
\def\pV{\mathbf{V}}

\def\xict{\xi^{\,c,t}}

\def\M{M}
\def\ur{u^r}
\def\uh{u^h}

\def\Nmd{N_{\rm model}}

\def\mf{\mathbf{f}}

\def\mA{\mathbf{A}}
\def\mQ{\mathbf{Q}}

\def\mS{\mathbf{S}}

\def\lmut{\tilde{\mathbf{u}}_{\xi_t}}

\def\lmek{\tau_{\xi^{(j\,)}_t}}

\def\lano{\gamma_{t}}
\def\ltano{\tilde{\gamma}_{t}}

\def\lspan{\textrm{span}}
\def\Qe{Q_{\emptyset}}
\def\POD{\textrm{POD}}
\def\covl{\alpha}
\def\QI{Q^0}

\def\xitj{\xi_t^{(j)}}

\def\Rt{\color{black}}
\def\bolde{\epsilon}

\def\md{{\mathbf d}}

\begin{document}

\begin{frontmatter}
\title{An adaptive reduced basis ANOVA method for high-dimensional 
Bayesian inverse problems}

\author[QL]{Qifeng Liao}
\ead{liaoqf@shanghaitech.edu.cn}
\address[QL]{School of Information Science and Technology, ShanghaiTech University, Shanghai 201210, China}

\author[JL]{Jinglai Li\corref{mycorrespondingauthor}}
\cortext[mycorrespondingauthor]{Corresponding author}
\ead{Jinglai.Li@liverpool.ac.uk}
\address[JL]{Department of Mathematical Sciences, University of Liverpool, Liverpool L69 7XL, UK}

\begin{abstract}
In Bayesian inverse problems sampling the posterior distribution is often a challenging task when the underlying models are computationally intensive. 
To this end, surrogates or reduced models are often used to accelerate the computation. However, in many practical problems, the parameter of interest
can be of high dimensionality, which renders standard model reduction techniques infeasible. 
In this paper, we present an approach that employs the ANOVA decomposition method to reduce the model with respect to the unknown parameters, 
and the reduced basis method to reduce the model with respect to the physical parameters. Moreover, we provide an adaptive scheme within the MCMC iterations,
to perform the ANOVA decomposition with respect to the posterior distribution. With numerical examples, we demonstrate that 
the proposed model reduction method can significantly reduce the computational cost of Bayesian inverse problems,
without sacrificing much accuracy.   
\end{abstract}


\end{frontmatter}


\pagestyle{myheadings}
\thispagestyle{plain}

\section{Introduction} \label{section_intro}
Inverse problems arise from many fields of science and engineering---whenever parameters of interest must be estimated from indirect observations~\cite{tarantola2005inverse}.
The Bayesian inference method has become increasingly popular as a tool to solve inverse problems~\cite{tarantola2006popper,kaipio2006statistical}. The popularity of the method is largely due to its ability to quantify the uncertainty in the solution obtained.  Simply speaking, in the Bayesian framework, the parameters of interest are cast as random variables to which a prior distribution is assigned, and then the posterior distribution of the parameters conditional on the observed data is computed via the Bayes' rule. 
The posterior distribution thus provides a probabilistic characterization of the parameters of interest in such problems. 

Though the idea behind the Bayesian inference method is quite straightforward, the computation of the posterior distribution often poses challenges. 
In most practical problems, the posterior distributions do not admit a closed-form expression and must be computed numerically. To this end, the Markov chain Monte Carlo (MCMC) method~\cite{robert2005monte} is often used to compute the posterior distributions. 
In particular, the MCMC method draws samples from the posterior distribution and then any posterior statistics can be evaluated with the obtained samples.  
As will be shown later, the MCMC method requires to repeatedly evaluate the likelihood function and each evaluation involves a full simulation of the forward function, i.e., the mapping from the parameter of interest to the observables.
In many practical problems,  such as the seismic inversion~\cite{virieux2009overview} and the inverse groundwater modeling~\cite{yeh1986review},  
the forward functions are often described by computationally intensive partial differential equations (PDEs). 
On the other hand, often a rather large number of samples are required to accurately estimate certain posterior moments. 
In this case, the total cost of the MCMC simulation may become prohibitively high. 

To improve the overall computational efficiency of the MCMC simulation, we can reduce either the number of samples required or  the cost for generating each sample. 
The first option is essentially to develop more efficient sampling schemes, which is not in the scope of  this work. 
We here consider the second option, i.e., to reduce the cost for computing each sample.
To this end, a natural idea is to construct computationally inexpensive surrogates/reduced models and use 
them in the MCMC procedure. Substantial efforts have been made toward this direction and various types of surrogate models have been 
used to approximate the forward functions, most notably, the polynomial chaos expansion (PCE)~\cite{marzouk2009stochastic,marzouk2009dimensionality,marzouk2007stochastic,nagel2016spectral,yan2015stochastic},
the Gaussian process surrogates~\cite{bilionis2013solution,kennedy2001bayesian,wang2017adaptive}, the sparse grid interpolation~\cite{ma2009efficient}, and the reduced order models (ROM)~\cite{cui2015data,galbally2010non,lieberman2010parameter,wang2005using,CMCS-CHAPTER-2008-001}. 
The performances of these methods (especially the PCE and the ROM) for accelerating the Bayesian computation is detailedly compared 
and discussed in \cite{frangos2010}.  

The surrogate or reduced model based methods have been successfully applied to a large variety of inverse problems, resulting in significant computational saving of the Bayesian inference.
Despite the success, the applicability of these methods is often ultimately limited by the dimensionality of the unknown parameters. 
In many real-world applications, the unknown parameters are often of very high dimensionality: for example, in groundwater modelling one may want to estimate the hydraulic conductivity,
and in seismic inversion it is the wave velocity that one is interested in;  
in these problems, the unknowns are spatial fields, and if we represent the unknown fields with mesh grid points, the resulting inverse problems can be of tens of thousands or more dimensions. 
Doing Bayesian inference directly for such problems is often not possible, 
and in practice it often requires dimension reduction for the input space to make the inference feasible.  
In particular,  the truncated Karhunen-Loeve (KL) expansion is often used to represent the unknown field that we want to estimate~\cite{marzouk2009dimensionality,li2015note} to reduce the dimensionality.
However, in many practical problems,  the unknown fields are rough, and in this case one still needs to use a rather large number 
of KL modes to represent it. 
Constructing surrogate or reduced models for such high-dimensional problems directly is a  rather challenging task. 

The main purpose of this work is to provide an approach to tackle the dimensionality issue and construct reduced models for such problems. 
Specifically, we focus on the Analysis of Variance (ANOVA) methods~\cite{fisher,sob03,caogun09,wintar09,gaohes10,zhachokar11}. 
The ANOVA methods, 
which are proposed for efficiently solving high-dimensional forward uncertainty quantification (UQ) 
problems, aim to decompose a high-dimensional parameter space into a union of low-dimensional spaces, such that standard 
surrogate/reduced modelling strategies can be applied. For example, these include ANOVA based stochastic collocation 
\cite{maza10, yanlin11}, ANOVA multi-element collocation \cite{fookar09}, 
and reduced  basis ANOVA \cite{zhahes16,liaolin16,choelman18}.  
However,
how to develop an efficient ANOVA approach for high-dimensional Bayesian inversion still remains an open question. 
The main difficulty here is that, conducting ANOVA decomposition of high-dimensional models 
requires the knowledge of the distribution of the input parameters, 
which in the Bayesian inverse problems is the posterior that we want to compute. 
An approximate solution is to perform ANOVA decomposition with respect to the prior distribution, 
but the prior based ANOVA decomposition is often inefficient,
especially when the prior is significantly different from the posterior. 
Thus, we develop an adaptive  reduced basis ANOVA (RB-ANOVA) algorithm which allows us to construct a reduced model with respect to the posterior 
distribution, which, as is illustrated by numerical examples, is more efficient than that constructed based on the prior. 
To summarize, the main contributions of this work are two-fold: first we propose to use the RB-ANOVA model to accelerate the MCMC 
simulations for high-dimensional Bayesian inverse problems; 
second, we develop an adaptive scheme to construct the RB-ANOVA model with respect to the posterior distribution.

The rest of the paper is organized as follows. In Section~\ref{section_inference} we describe the formulation of the Bayesian inverse problems that will be considered in this work. 
In Section~\ref{section_rbANOVA} we provide a scheme for constructing the RB-ANOVA model,  and in Section~\ref{section_rbmcmc} we present our new RB-ANOVA based Markov chain Monte Carlo (RB-ANOVA-MCMC) algorithm,
which adaptively constructs the RB-ANOVA model with respect to the posterior distribution within the MCMC iterations. 
In Section~\ref{section_numerical}, with numerical experiments we demonstrate that the proposed adaptive RB-ANOVA method can significantly accelerate the Bayesian computation. 
Finally some concluding remarks are offered in Section~\ref{section_conclude}.

\section{Bayesian inverse problems}\label{section_inference}
In this section we describe the problem setup that is used in this work.
Suppose that we are interested in an $M$-dimensional parameter 
$\xi=[\xi_1,\ldots,\xi_\M]^T \in \dsR^{M}$,  
and we want to estimate it from some observed data $\-d$. 
Moreover we assume that there exists a forward model $\-G$ that maps the unknown parameter $\xi$ to the data $d$:
 \begin{eqnarray}\label{eq:setting}
 d =\-G(\xi)+\bolde,
\end{eqnarray}
 where $\bolde$ is the measurement noise. 
 Let $\pi_\epsilon(\bolde)$ be the distribution of $\bolde$, and one can obtain the distribution of $d$ conditional on $\xi$:
\begin{equation}
\pi(\-d|\xi) = \pi_\epsilon(\-d-\-G(\xi))\,.\label{e:lh}
\end{equation}
In a Bayesian formulation, one assigns a prior distribution $\pi(\xi)$ on $\xi$ encoding the prior knowledge 
on the parameter of interest, and the posterior $\pi(\xi|\-d)$ can then be calculated via Bayes' rule:
\begin{equation}
\pi(\xi|\-d) = \frac{\pi(\-d|\xi)\pi(\xi)}{\int \pi(\-d|\xi)\pi(\xi) d\xi\,},
\label{e:bayesrule}
\end{equation}
where the denominator is a normalization constant that makes the posterior a well-defined probability distribution. 
We note here that, in practice it is usually reasonable to assume that the sought parameters 
are in a (sufficiently large) bounded region, and thus in this paper, we shall restrict our attention to the situation that  
the prior $\pi(\xi)$ has a bounded and connected support.
Without loss of generality, we then assume the support of $\pi(\xi)$ is  $I^\M$ where $I:=[-1,1]$ throughout this work.
 
As is mentioned earlier, one frequently employs the MCMC simulation to sample the posterior distribution.
Simply speaking, the MCMC method constructs a Markov chain which asymptotically converges to the posterior distribution.
In this work, we adopt the popular 
Metropolis-Hastings (MH) MCMC algorithm outlined in Algorithm~\ref{alg:mh},
to generate  $N$ samples $\{\xi^{(1)},\ldots,\xi^{(N)}\}$ of the posterior of $\xi$.
In Algorithm~\ref{alg:mh},  $\pi(\cdot |\-\xi^{(j)})$ on line $3$ is a given proposal distribution
which may be a multivariate normal distribution with mean $\xi^{(j)}$, and $U[0,1]$ on line 5 refers to
the uniform distribution on $[0,1]$.
\begin{algorithm}
\caption{The standard MH algorithm}
\label{alg:mh}
\begin{algorithmic}[1]
\STATE  Initialize the chain at $\xi^{(1)}$.
\FOR{$j=1:N-1$}
                  \STATE Draw $\-\xi^*\sim \pi(\cdot |\-\xi^{(j)})$.                  
                  \STATE Compute the acceptance ratio 
$a = 
\min\left(1, \frac{\pi_\epsilon\left(d-G_{\pP}\left(\xi^*\right)\right)\pi(\xi^*)}{\pi_\epsilon\left(d-G_{\pP}\left(\xi^{(j)}\right)\right)\pi(\xi^{(j)})}
\frac{\pi\left(\xi^{(j)}|\xi^*\right)}{\pi\left(\xi^{*}|\xi^{(j)}\right)}\right).$

                  \STATE Draw $\rho\sim U[0,1]$.
                  \IF{$ \rho<a$} 
  \STATE                  Let $\-\xi^{(j+1)} = \-\xi^*$, 
                  \ELSE 
                  \STATE Let $\-\xi^{(j+1)}=\-\xi^{(j)}$.
                  \ENDIF
\ENDFOR

\end{algorithmic}
\end{algorithm}

It can be seen from the algorithm that, each MCMC iteration requires 
an evaluation of the computationally expensive forward function $G(\cdot)$
(on line $4$ of Algorithm~\ref{alg:mh}), 
which renders the MCMC procedure formidably expensive. 
In what follows we provide a reduced basis ANOVA based method to accelerate the MCMC computation.

\section{The RB-ANOVA  method}\label{section_rbANOVA}
To begin with, details of the forward model considered in this paper are addressed as follows.
Let $D$ denote a spatial domain (a subset of $\dsR^2$ or $\dsR^3$) which is  bounded, connected and with a
polygonal boundary $\partial D$,
and $x\in D$ denote a spatial variable.
The physics of problems considered 
are governed by 
a PDE over the spatial domain $D$ and
boundary conditions on the boundary $\partial D$,
which are stated as: find
$u(x,\xi)$ mapping $D\times I^M$ to $\dsR$, such that
\begin{subequations} \label{e:pdemodel}
\begin{align}
\pL\left(x,\xi;u\left(x,\xi\right)\right)=f(x)  \qquad
&\forall \left(x,\xi\right) \in D\times I^M,
\label{spdexi1}\\
\pB\left(x,\xi;u\left(x,\xi\right)\right)=g(x) \qquad
&\forall \left(x,\xi\right)\in \partial D\times I^M,
\label{spdexi2}
\end{align}
\end{subequations}
where $\pL$ is a partial differential operator and $\pB$ is a boundary
operator, both of which can depend on the unknown parameter $\xi$. Here $f$ is the source function and $g$ 
specifies the boundary conditions. Through specifying an observation operator $\pC$, e.g., 
taking solution values at several grid points, 
we write the overall forward model as $G(\xi):=\pC(u(x,\xi))$. 
It is clear that each evaluation of the forward function requires to solve the PDE~\eqref{e:pdemodel}, 
and this procedure needs to be performed repeatedly 
in the MCMC iterations. 
As discussed earlier, we shall construct computationally inexpensive reduced models and use 
them in the MCMC iteration to accelerate the computation.  
However,  when the parameter of interest is high-dimensional, 
constructing reduced models are rather challenging. 
In this work, the ANOVA decomposition approach is used to decompose the model so that the reduced model
construction becomes feasible.    
The construction of the RB-ANOVA surrogate for the forward models 
is discussed in this section, which is an extension of the procedure outlined in \cite{liaolin16}, 
and the application of it to Bayesian inversion is presented in the next section.

\subsection{ANOVA decomposition}\label{sec:anova}
We present the ANOVA decomposition method in a generic setting.
Namely, suppose that we have a computationally intensive function $u(\-x,\xi)$ 
where $x\in D$ is the physical variable 
and $\xi\in I^M$ is the random variable, 
and the goal here is to construct a reduced model (or approximation) of $u(x,\xi)$ with respect to the random variable $\xi$. 

To proceed, the notation for indices are first set up following \cite{gaohes10,liaolin16}. 
In general, any subset of $\{1,\ldots,M\}$ denotes an index. 
For an index $t\subseteq \{1,\ldots,M\}$, $|t|$ denotes the cardinality of $t$,
and we define $|t|=0$ for $t=\emptyset$. 
For an index $t\neq \emptyset$, we sort its elements in ascending order and
express it as $t=(t_{1},\ldots,t_{|t|})$ with $t_{1}<t_{2}\ldots <t_{|t|}$.
In addition, we also call $|t|$ the (ANOVA) order of $t$, and call $t$ a $|t|$-th order index.
For a given ANOVA order $i=0,\ldots,M$, the following index sets are defined
\begin{eqnarray*}
\pd_i&:=&\left\{t\,|\> t\subset\{1,\ldots,M\}, \>|t|=i \right\},\\
\pd^{\star}_i&:=&\cup_{j=0,1,\cdots,\,i}\pd_j,\\
\pd&:=&\pd^{\star}_M=\cup_{j=0,1,\cdots,\,M}\pd_j.
\end{eqnarray*}
The sizes of the above sets (numbers of elements that they contain) are denoted by 
$|\pd_i|$, $|\pd^{\star}_i|$ and $|\pd|$ respectively. From the above definition, $\pd_0=\{\emptyset\}$ and $|\pd_0|=1$. 
For a given index $t=(t_1,\ldots, t_{|t|})\in \pd$ with $|t|>0$, 
$\xi_t$ denotes a random vector collecting components of $\xi$ associated with
$t$, i.e., $\xi_t:=[\xi_{t_1},\ldots,\xi_{t_{|t|}}]^T \in I^{|t|}$, and we denote
the (marginal) prior probability density function of $\xi_t$ by $\pi_t(\xi_t)$ and its (marginal) posterior
probability density function by $\pi_t^*(\xi_t):=\pi_t(\xi_t| d)$.

While noting that there are other strategies to implement the ANOVA decomposition~\cite{maza10,yanlin11,zhachokar12}, 
here we  adopt the so-called anchored ANOVA method following \cite{maza10,yanlin11,zhachokar11,liaolin16}.
In this method, one first selects an anchor point $c=[c_1,\ldots,c_M]^T\in I^M$, and then decomposes
the function $u(x,\xi)$ with respect to $\xi = [\xi_1,\dots,, \xi_M]^T\in I^M$ as, 
\begin{eqnarray}
u(x,\xi)&=& u_{0}(x)+u_1(x,\xi_1)+\ldots+u_{1,2}(x,\xi_{1,2})+\ldots \nonumber \\
&=&\sum_{t\in \pd} u_t(x,\xi_t), \label{anova}
\end{eqnarray}
where we denote $u_{\emptyset}(x,\xi_{\emptyset}):=u_0(x)$ for convenience, and each term 
in \eqref{anova} is specified as 
\begin{subequations} \label{e:ut}
\begin{eqnarray}
u_{\emptyset}(x,\xi_{\emptyset})&:=&u_0(x):=u(x,c),\label{anova_0}\\
u_t(x,\xi_t)&:=&u(x,c,\xi_{t})-\sum_{s\subset t}u_s(x,\xi_s).\label{anova_t}
\end{eqnarray}
\end{subequations}
In the equation above, we have $\xi_\emptyset=c$, and $u(x,c,\xi_t)$ is defined as, 
\begin{eqnarray}
u(x,c,\xi_{t})&:=&u\left(x,\xict\right),\nonumber
\end{eqnarray}
where 
\begin{subequations}
\label{xict}
\begin{eqnarray}
\xict&:=&[\xict_1,\ldots,\xict_M]^T,\\
\xict_i&:=&\left\{\begin{array}{ll} 
                     c_i & \textrm{for } i\in \{1,\ldots,M\}\setminus t\\
                     \xi_i & \textrm{for } i \in t 
                     \end{array}\right. .
\end{eqnarray}
\end{subequations}
In what follows, $u_t(x,\xi_t)$ is called a \emph{child term} of $u_s(x,\xi_s)$ if $s\subset t$. 
It should be clear that the decomosition~\eqref{anova} is exact and so itself does not provide us
 a reduced model of the solution $u(x,\xi)$. 
However, as  
discussed in \cite{maza10,yanlin11,liaolin16}, an efficient reduced model can be obtained 
if one only keeps a small number of \emph{active} terms in \eqref{anova}.
We will discuss how to select the active terms later. 
For now supposing that we have selected the active terms,
the sets consisting of selected important indices at each order
are denoted by $\pP_i \subseteq \pd_i$ for $i=0,\ldots,M$. 
We then define $\pP^{\star}_i:=\cup_{j=0,\ldots,i}\pP_j$ and $\pP:=\pP^{\star}_M$.
A reduced model of the solution $u(x,\xi)$ is obtained:  
\begin{eqnarray}
u\left(x,\xi\right)\approx u_{\pP}\left(x,\xi\right):= \sum_{t\in \pP} u_t\left(x,\xi_t\right), \label{anova_adaptive}
\end{eqnarray}
where $u_t$ is defined in (\ref{anova_t}).
In the following, $u_{\pP}\left(x,\xi\right)$
is called the ANOVA model (or approximation) of $u(x,\xi)$. 

For selecting the active terms (or indices) in the ANOVA model, 
the prior distribution $\pi_t$ of $\xi_t$ is given in advance in this Bayesian inference setting,
and thus a natural idea is to construct the selection criterion using some prior statistics.
While they are not optimal choices, the prior statistics are used to illustrate 
the methods in this section, and optimal selection criteria based on posterior distributions 
are presented in our new algorithm in the next section.  
To this end, we adopt the relative mean approach used in \cite{liaolin16}, 
while noting that other choices are also possible~\cite{maza10,yanlin11}. 
Specifically, recalling that the prior mean of $u_t$ is
\[\pE\left(u_t\right):=\int_{I^{|t|}} u_t\left(x,\xi_t\right) \pi_t\left(\xi_t\right)\,d\xi_t,\]
we define the relative mean value to be
\[
\lano:=\frac{\left\|{\pE}(u_t)\right\|_{0,D}}{\left\|\sum_{s\in \pP^{\star}_{|t|-1}}{\pE}\left(u_s\right)\right\|_{0,D}},
\]
where $\|\cdot\|_{0,D}$ denotes the $L^2$ function norm over region $D$. 
In practice, the prior expectation $\pE\left(u_t\right)$ can be computed with a Monte Carlo (MC) estimator:
\begin{eqnarray}
\ppE \left(u_t\right):=\frac{1}{N}\sum_{j=1}^Nu_t\left(x,\xitj\right), \label{mean_mc}
\end{eqnarray}
where $ \{\xitj\}_{j=1}^N$ are $N$ samples drawn from $\pi_t$,
and as a result,  the relative mean value $\lano$ can be approximated by
\begin{eqnarray}
\tilde{\lano}:=\frac{\left\|\tilde{\pE}(u_t)\right\|_{0,D}}{\left\|\sum_{s\in \pP^{\star}_{|t|-1}}\tilde{\pE}\left(u_s\right)\right\|_{0,D}}. 
\label{etamc}
\end{eqnarray}
Here we call a term $u_t$ \emph{important} if the associated relative mean estimate is larger than a prescribed threshold value
$tol_{anova}$. 
The set of active terms at each order is selected with the following procedure. 
Namely, suppose that $\pP_{i}$ is given,  
and one first selects all  important terms at order $i$, yielding the index set \[\tilde{\pP}_{i}:=\{t\,| \,t\in\pP_{i} \textrm{ and }\lano\geq tol_{anova} \},\]
which is  a subset of $\pP_{i}$. After that, as discussed in \cite{maza10}, the index set at order $i+1$ is constructed by 
\begin{eqnarray}
\pP_{i+1}:=\left\{t\,| \,t\in\pd_{i},\textrm{ and any } s\subset t \textrm{ with $|s|=i$ satisfies } s\in \tilde{\pP}_i\right\}.\label{pkk}
\end{eqnarray}
That is, 
if a term is found unimportant, the term itself is not removed from the ANOVA model, 
but all its child terms are removed for the next order. 
To start the procedure, we set $\tilde{\pP}_0=\pP_0=\pd_0=\emptyset$. 
On the other hand, the procedure terminates automatically if no active term is found for the next order. 
The studies in \cite{maza10,yanlin11} indicate that for most realistic physical systems the size of $\pP$ is usually much smaller than that of $\pd$, and moreover, $\pP$ 
may only contain low order terms.  

\subsection{The RB  approximation}\label{section_rb}
In the present problem, $u(x,\xi)$ is the solution of the parameterized equation~\eqref{e:pdemodel}. 
As mentioned in the previous section, 
the ANOVA decomposition method yields a reduced model in the random parameter space.
Here we 
discuss how to perform model reduction with respect to the physical parameter $x$, with the reduced basis (RB) method.   

First, to use the ANOVA model~\eqref{anova_adaptive}, the terms $u(x,c,\xi_t)$ in 
\eqref{e:ut} for all $t\in\pP$ need to be computed.   
Here, $u(x,c,\xi_{t})$ is the solution of the following equations:
 \begin{subequations} \label{e:local}
\begin{align}
\pL_t\left(x,\xi_t;u\left(x,c,\xi_t\right)\right)=f(x)  \qquad 
&\forall \left(x,\xi_t\right) \in D\times I^{|t|},
\label{spdexi1_t}\\ 
\pB_t\left(x,\xi_t;u\left(x,c,\xi_t\right)\right)=g(x) \qquad  
&\forall \left(x,\xi_t\right)\in \partial D\times I^{|t|},
\label{spdexi2_t}
\end{align}
\end{subequations}
where $u(x,c,\xi_t)$ is defined by \eqref{xict} and $\pL_t$ and $\pB_t$ are defined through putting \eqref{xict} into \eqref{e:pdemodel}. 
Eqs.~\eqref{e:local} are referred to as a (parametrically) $|t|$-dimensional {\em local} problem,
while the global problem is Eqs.~\eqref{e:pdemodel}.
It is easy to see that, if $u(x,c,\xi_t)$ is evaluated by directly solving 
the local problem~\eqref{e:local} with the same strategy for solving \eqref{e:pdemodel}, 
evaluating the ANOVA model~\eqref{anova_adaptive} is actually much more expensive than solving the global problem~\eqref{e:pdemodel} directly.
This is because that the ANOVA model requires to solve the local problem multiple times and a full solve of the local problem is about as costly as that of the global problem.  
Thus, to make the ANOVA model useful for our problem, a reduced model for the local 
problem~\eqref{e:local} needs to be constructed, so that it can be solved more efficiently. 
We construct such a  model using the RB method. 

We start with the finite element approximation of the local problem~\eqref{e:local}.
In general, the variational form of the deterministic problem \eqref{e:local} 
corresponding to a given realization of $\xi_t$ is given by
$\ssl_{\xi_t}(u(x,c,\xi_t),v)=\ssr(v)$.
Given a finite element space $X^h$ with $N_h$ degrees of freedom, 
a finite element formulation seeks a solution $u^h(x,c,\xi_t)\in X^h$ such that
\begin{eqnarray}
\ssl_{\xi_t}\left(u^h(x,c,\xi_t),v\right)=\ssr(v),\quad \forall v\in X^h. \label{fem}
\end{eqnarray}
As usual, a finite element solution $u^h$ is referred to as a {\em snapshot}. 
Next, the reduced basis (RB) approximation is stated as: given a set of reduced basis functions
$Q_t:=\{q_t^{(1)},\cdots,q^{(N_r)}_t\}\subset X^h$, 
find $u^r(x,c,\xi_t)\in \lspan\{Q_t\}$ such that 
\begin{eqnarray}
\ssl_{\xi_t}\left(u^r(x,c,\xi_t),v\right)=\ssr(v),\quad \forall v\in \lspan\{Q_t\}. \label{rb}
\end{eqnarray}

Two standard methods are used to generate the reduced bases $Q_t$
for all $t\in\pP$ in this paper.
The first one is the proper orthogonal decomposition (POD) \cite{sirovich87,holmes96,gunpetsha07},
which can be briefly reviewed as follows.
For a given finite sample set $\Xi \subset I^{|t|}$ with size $|\Xi|$, a finite snapshot set is defined by
\begin{eqnarray}
S^t_{\Xi}:=\left\{u^h\left(x,c,\xi_t\right),\,\xi_t\in \Xi \right\}.
\end{eqnarray}
The matrix form of $S^t_{\Xi}$ is denoted by 
$\mS^t_{\Xi}\in \dsR^{N_h\times |\Xi|}$,
i.e., each column of $\mS^t_{\Xi}$ is the vector of basis function coefficients 
of a finite element solution. 
Assuming $|\Xi|<N_h$,
let $\mS^t_{\Xi}=U \Sigma V^T$ denote the 
singular value decomposition (SVD) of $\mS^t_{\Xi}$, where $U=(\mathbf{q}_1,\cdots,\mathbf{q}_{|\Xi|})$ and $\Sigma={\rm diag}
(\sigma_1,\cdots,\sigma_{|\Xi|})$ with $\sigma_1\geq \sigma_2\geq \cdots \geq \sigma_{|\Xi|}\geq 0$. 
The basis $Q_t$ is then given by the first $k$ left singular vectors $(\mathbf{q}_1,\ldots,\mathbf{q}_k)$, of which the corresponding singular 
values are greater than some given tolerance $tol_{pod}$, 
i.e., $\sigma_{k}/\sigma_{1}> tol_{pod}$ but $\sigma_{k+1}/\sigma_{1}\leq tol_{pod}$. 
As usual, to simplify the later presentation, 
this POD procedure for generating  $Q_t$ through $S^t_{\Xi}$ is denoted by $Q_t:=\POD(S^t_{\Xi})$. 

The second one is the greedy sampling method
\cite{verpat02,nguver05,haaohl08,buiwillcox08,boybri10,patrozbook,quamanbook}.
This method is to adaptively select parameter samples, where errors between the reduced approximation and
the finite element approximation are large. 
To assess the errors, we use the residual error indicator which is also adopted by 
\cite{elmanliao,liaolin16,powellnewsum16,zhangwebster17}.
Following our notation in \cite{elmanliao},
when considering linear 
PDEs, the algebraic system associated with (\ref{fem}) can be written as 
$\mA_{\xi_t}\mathbf{u}_{\xi_t}=\mathbf{f} $ where 
$\mA_{\xi_t}\in\dsR^{N_h\times N_h}$, and $\mathbf{u}_{\xi_t},\,\mathbf{f} \in \dsR^{N_h}$. 
The algebraic system of the reduced basis approximation (\ref{rb}) 
can be written as $\mathbf{Q}_t^T\mA_{\xi_t}\mathbf{Q}_t\lmut=\mathbf{Q}_t^T\mathbf{f},
\label{matrix_rb}$
where $\lmut\in \dsR^{N_r}$  gives a reduced basis solution and 
$\Rt \mQ_t\in \dsR^{N_h\times N_r}$ is the matrix form of the reduced basis $\Rt Q_t=\{q_1,\ldots,q_{N_r}\}$, 
i.e., each column of $\mQ_t$ is the vector of nodal coefficient values associated 
with each $q_i$, $i=1,\ldots,N_r$. The residual indicator is defined by
\begin{eqnarray}
\tau_{\xi_t}:=\frac{\|\mA_{\xi_t}\mQ_t\lmut-\mf\|_2}{\|\mf\|_2}. \label{res}
\end{eqnarray} 
With this residual indicator, the greedy sampling procedure can be stated as follows.
First, take the first sample $\xi_t^{(1)}$ from a given sample set $\Xi$ and initialize the reduced basis as $Q_t:=\{u^h(x,c,\xi^{(1)}_t)\}$.
Second, for each $\xi_t\in \Xi$, compute the residual error indicator $\tau_{\xi_t}$
using the current reduced basis $Q_t$, and if $\tau_{\xi_t}$ is larger than some given tolerance, compute 
the snapshot $u^h(x,c,\xi_t)$ and augment $Q_t$ with $u^h(x,c,\xi_t)$.
The second step is repeated until  $N_r$ snapshots are obtained.

\subsection{The RB-ANOVA model}\label{sec:rbanova}
With the local problem~\eqref{e:local}  solved by the RB method, we obtain a  RB-ANOVA model:
\begin{eqnarray}
\ur_{\pP}\left(x,\xi\right):= \sum_{t\in \pP} \ur_t\left(x,\xi_t\right), \label{rb_anova}
\end{eqnarray}
where
\begin{subequations}
\label{rb_anova_0t}
\begin{eqnarray}
\ur_{\emptyset}(x,\xi_{\emptyset})&:=&u^h(x,c),\label{rb_anova_0}\\
\ur_t(x,\xi_t)&:=&\ur(x,c,\xi_{t})-\sum_{s\subset t}\ur_s(x,\xi_s).\label{rb_anova_t}
\end{eqnarray}
\end{subequations}
In \eqref{rb_anova_0t}, $u^r(x,c,\xi_t)$ is the RB solution of the local problem~\eqref{rb}, and $u^h(x,c)$ is 
the snapshot at the anchor point (i.e., the solution of \eqref{fem} with $t=\emptyset$). 
Constructing the RB-ANOVA model in our setting is equivalent to generating four pieces of data:
the anchor point $c$, the snapshot $u^h(x,c)$ at the anchor point, 
the index set $\pP$, and the reduced basis $Q_t$ for each $t\in \pP$. 
We call these data the RB-ANOVA model data. 
With them, a  RB-ANOVA approximation $\ur_{\pP}(x,\xi)$ at any input sample point $\xi\in I^M$ can be cheaply computed.
The procedures for generating the RB-ANOVA data $\{c,\, \uh(x,c),\, \pP,\, \{Q_t\}_{t\in\pP}\}$ 
are as follows.

First, suppose that we are given a set of  realizations of the random variable $\xi$, denoted by $\Xi$.
As discussed in \cite{gaohes10}, for a given distribution of $\xi$,
the optimal anchor point $c$ with respect to this distribution is its mean point.
However, the goal of this work is to generate samples for the posterior distribution, 
of which the exact mean point is not admitted. As an alternative, the anchor point in this work 
is taken to be the sample mean of $\Xi$.

We set $\pP_0:=\{\emptyset\}$, and compute the snapshot  $u^h(x,c)$. 
The zeroth order RB is constructed using this snapshot $Q_{\emptyset}:=\{u^h(x,c)\}$,
and the mean estimate for the zeroth order ANOVA term is set to  
$\ppE(u_{\emptyset}):=u^h(x,c)$.
Moreover, it is easy to see that $\tilde{\pP}_0={\pP}_0$, which immediately implies that $\pP_1:=\{1,\ldots,M\}$. 
Now we consider an ANOVA order $i\geq 1$.
That is, given the index set $\pP_i$ and the reduced bases for order $i-1$, $\{Q_s\}_{s\in \pP_{i-1}}$, we need to find  the set $\tilde{\pP}_{i}$ and the reduced bases $\{Q_t\}_{t\in \pP_{i}}$. 
Now recall that, the set $\tilde{\pP}_i$ is obtained by estimating the relative means with MC approximation. 
It should be clear that here if the Monte Carlo samples of $\{u_t(x,\xi_t^{(j)})\}_{j=1}^N$ for each $t\in \pP_i$ are computed 
with the PDE model with finite elements, the total cost may become prohibitively high.  
To reduce the cost, we consider the reduced basis MC method which incorporates 
 greedy RB methods in MC simulations \cite{boybri10}, and extend it to yield both the set $\tilde{\pP}_i$ 
and the reduced bases $\{Q_t\}_{t\in \pP_{i}}$ with low costs.

To start the greedy procedure, the hierarchical approach introduced in \cite{liaolin16} is used 
to initialize the reduced basis $Q_t$ for $t\in \pP_i$, which reuses the bases generated at the previous order
based on  the nested structure of ANOVA indices:
\begin{enumerate}
\item  grouping all reduced basis functions associated with subindices of $t$ with
order $|t|-1$ together, we define $\QI_t:=\cup_{s\in \Lambda_t} Q_s$ 
where $\Lambda_t:=\{s\,|\,s\in\pP_{|t|-1}  \textrm{ and } s\subset t\}$;
\item  we apply POD to $\QI_t$ to result in an orthogonal basis to serve as an initialization of $Q_t$, 
i.e., we initially set $Q_t:=\POD(\QI_t)$ (details of POD are discussed in Section \ref{section_rb}).
\end{enumerate}
After the initial basis is generated, a sample set of $\xi_t$ for $t\in \pP_i$ needs to be specified to 
conduct the MC simulation. Since the sample set $\Xi$ is given for the global parameter $\xi$
and each $\xi_t$ for $t\in\pP$ is a collection of components of $\xi$, it is trivial to define a sample set of
$\xi_t$ by a collection of the components of samples in $\Xi$, i.e., the samples of $\xi_t$ are taken to 
be $\Xi_t:=\{\xi^{(j)}_t,\xi^{(j)}\in\Xi \textrm{ for } j=1,\ldots, |\Xi| \}\subset I^{|t|}$.
Then, looping over the sample points, we compute the reduced solution $u^r(x,c,\xi^{(j)}_t)$ (see (\ref{rb}))
for each $\xi^{(j)}_t\in \Xi_t$,
and the residual indicator $\lmek$ (see (\ref{res})):
\begin{enumerate}
\item if the residual indicator is smaller than a given tolerance $tol_{rb}$, use $u^r(x,c,\xi^{(j)}_t)$  to serve as 
a MC solution sample;
\item if the residual indicator is larger than or equal to $tol_{rb}$, compute the snapshot $u^h(x,c,\xi_t^{(j)})$
through solving (\ref{fem}), use the snapshot to
serve as a MC solution sample and update the reduced basis $Q_t$ with this snapshot.
\end{enumerate}
When all $|\Xi_t|$ MC samples are generated through the above greedy approach, we 
compute the relative mean values using \eqref{etamc} and construct the important index set $\tilde{\pP}_i$, which consequently yields $\pP_{i+1}$.   
As is mentioned in Section~\ref{sec:anova}, the above procedure is repeated until  $\pP_{i+1}= \emptyset$.
This RB-ANOVA procedure  is 
formally stated in Algorithm \ref{alg_rbanova}. 
It should be noted that this algorithm only requires a set of realizations of $\xi$, $\Xi$, as its input, 
and this is an important property for the adaptive algorithm that will be presented in the next section. 
We also note that, a major difference between Algorithm~\ref{alg_rbanova} and that developed in \cite{liaolin16} is that, in \cite{liaolin16} the RB-ANOVA model
is constructed with the tensor grid collocation points, while here MC samples are used.

\begin{algorithm}
\caption{Constructing the RB-ANOVA model}
\label{alg_rbanova}
\begin{algorithmic}[1]       
\STATE Input: a finite sample set $\Xi:=\left\{\xi^{(j)},j=1,\ldots,|\Xi| \right\} \subset I^M $.
\STATE Compute the anchor point  $c:=\frac{1}{|\Xi|}\sum^{|\Xi|}_{j=1}\xi^{(j)}$.
\STATE Set $\pP_0:=\{\emptyset\}$, compute $u^h(x,c)$ (see \eqref{fem}) and set $u_{\emptyset}(x,\xi_{\emptyset}):=u^h(x,c)$.
\STATE Set $\Qe:=\{u^h(x,c)\}$, $\ppE(u_{\emptyset}):=u^h(x,c)$.
\STATE Set $\pP_1:=\{1,\ldots,M\}$, initialize $\pP:=\pP_0\cup\pP_1$, and let $i=1$.
\WHILE{$\pP_{i}\neq \emptyset$}
         \FOR{$t\in \pP_i$}
                 \STATE Construct $\QI_t:=\cup_{s\in \Lambda_t} Q_s$ where $\Lambda_t:=\{s\,|\,s\in\pP_{|t|-1}  \textrm{ and } s\subset t\}$.
                  \STATE
                   Initialize  $Q_t:=\POD\left(\QI_t\right)$, (see Section \ref{section_rb} for details of the POD method).
                   \STATE Construct the sample set $\Xi_t:=\left\{\xi^{(j)}_t,\xi^{(j)}\in\Xi \textrm{ for } j=1,\ldots, |\Xi|\right\}\subset I^i$.
                          \FOR{$j=1:\left|\Xi_{t}\right|$}
                                  \STATE Compute the reduced solution $u^r\left(x,c,\xi^{(j)}_t\right)$ through solving (\ref{rb}) and the error indicator $\lmek$ through (\ref{res}).
                                  \IF{$\lmek<tol_{rb}$}                          
                                      \STATE Set $u\left(x,c,\xi^{(j)}_t\right)=u^r\left(x,c,\xi^{(j)}_t\right)$ in \eqref{anova_t} 
                                      to obtain $u_t\left(x,\xi^{(j)}_t\right)$.
                                  \ELSE
                                     \STATE Compute the snapshot  $u^h\left(x,c,\xi^{(j)}_t\right)$  (see (\ref{fem})).
                                       \STATE Set $u\left(x,c,\xi^{(j)}_t\right)=u^h\left(x,c,\xi^{(j)}_t\right)$ in \eqref{anova_t} 
                                      to obtain $u_t\left(x,\xi^{(j)}_t\right)$.
                                     \STATE Augment the reduced basis $Q_t$ with $u^h\left(x,c,\xi^{(j)}_t\right)$,
                                     i.e.\ $Q_t=Q_t\cup\left\{u^h\left(x,c,\xi^{(j)}_t\right)\right\}$.
                                 \ENDIF
                  \ENDFOR
                  \STATE Compute $\ppE\left(u_t\right)$ using (\ref{mean_mc}) with samples $\left\{u_t\left(x,\xi_t\right),\,\xi_t\in\Xi_t\right\}$.
                  \STATE Compute the relative mean value
                  $\ltano={\left\|\ppE\left(u_t\right)\right\|_{0,D}}\,\left/\,{\left\|\sum_{s\in \pP^{\star}_{i-1}}\ppE\left(u_s\right)\right\|_{0,D}}\right.$.
         \ENDFOR
         \STATE Set $\tilde{\pP}_i:=\{t\,\left| \,t\in\pP_i, \textrm{ and }\ltano\geq tol_{anova} \right.\}$.
         \STATE Set $\pP_{i+1}:=\{t\,| \,t\in\pd_{i+1},\textrm{ and any } s\subset t \textrm{ satisfies } s\in \tilde{\pP}_i\}$.
         \STATE Update the index set $\pP:=\pP\cup\pP_{i+1}$ and update $i=i+1$.
\ENDWHILE
\STATE Output ANOVA model data: $\left\{c\,,\> u^h(x,c)\,,\> \pP \,,\> \{Q_t\}_{t\in\pP}\right\}$.
\end{algorithmic}
\end{algorithm}

We next discuss how to use the resulting RB-ANOVA model \eqref{rb_anova} 
to predict the system output $G(\xi)$ for an arbitrary input sample of $\xi$, 
as is required in the MCMC iterations. 
First,  we set $\ur_{\emptyset}(x,\xi_{\emptyset}):=\uh(x,c)$ as \eqref{rb_anova_0}.
Second, the reduced basis approximation $\ur(x,c,\xi_t)$
of the solution of each local system \eqref{e:local} for $t\in\pP$
is computed through solving (\ref{rb}) with the reduced basis $Q_t$. After that, $\ur_t(x,\xi_t)$ is computed 
through \eqref{rb_anova_t}, and the overall reduced basis ANOVA approximation $\ur_{\pP}(x,\xi)$
are computed through \eqref{rb_anova}. Finally, applying the given observation operator $\pC$ on $\ur_{\pP}(x,\xi)$,
the overall system output is estimated, i.e., we denote $G^r_{\pP}(\xi):=\pC(\ur_{\pP}(x,\xi))$.
This prediction procedure is summarized in Algorithm~\ref{alg_rbpre}.

\begin{algorithm}
\caption{Predication via reduced basis ANOVA model}
\label{alg_rbpre}
\begin{algorithmic}[1]
\STATE Input: a sample of $\xi$ and the RB-ANOVA model data $\left\{c,\,\uh(x,c),\, \pP,\, \{Q_t\}_{t\in\pP}\right\}$.
\STATE Set $\ur_{\emptyset}(x,\xi_{\emptyset}):=\uh(x,c)$.
\FOR{$t\in \pP\setminus\left\{\emptyset\right\}$}
                  \STATE Compute $\ur\left(x,c,\xi_t\right)$ through solving (\ref{rb}) with the reduced basis $Q_t$.
                  \STATE Obtain $\ur_t(x,\xi_t)$ through \eqref{rb_anova_t}.
\ENDFOR
\STATE Assemble $\ur_{\pP}\left(x,\xi\,\right)$ using \eqref{rb_anova}.
\STATE Set $G^r_{\pP}\left(\xi\, \right):=\pC\left(\ur_{\pP}\left(x,\xi\right)\right)$ where $\pC$ is the given observation operator.
\STATE Output: $G^r_{\pP}\left(\xi\,\right)$.
\end{algorithmic}
\end{algorithm}

\section{The adaptive RB-ANOVA method to accelerate MCMC}\label{section_rbmcmc}
In Section~\ref{sec:rbanova}, the schemes for constructing and using the RB-ANOVA surrogate for the forward  models
are presented. 
In the MCMC iterations, the computationally intensive finite element method can be replaced
with the RB-ANOVA model to reduce the computational cost. 
As discussed in Section~\ref{sec:rbanova}, a simple way of doing this is to construct the RB-ANOVA model with respect to the prior distribution before performing 
the MCMC simulation,  
which means that the sample set used to construct the reduced model in Algorithm \ref{alg_rbanova}  is generated from the prior distribution $\pi(\xi)$. 
An issue here is that,  the goal of the Bayesian inference is to sample according to the posterior distribution, 
and in this case, constructing the reduced model 
with respect to the prior distribution may become ineffective, especially for problems in which  
the posterior differs significantly from the prior~\cite{li2014adaptive}. 
Ideally one should construct the reduced model with respect to the posterior distribution for such problems, 
but this certainly can not be done in advance as the posterior is not available in advance.
To address the issue, we here present an algorithm that can adaptively construct the RB-ANOVA model according to the posterior distribution. 
Specifically, the new method updates the RB-ANOVA model inside the MCMC iterations, and 
for conciseness we shall refer to the whole procedure as the RB-ANOVA-MCMC alogorithm in the following.

In this section, the number of samples for generating the RB-ANOVA model is denoted by $\Nmd$, 
i.e., $|\Xi|=\Nmd$ on line 1 of Algorithm \ref{alg_rbanova}. 
To begin with, we construct an initial RB-ANOVA model using Algorithm~\ref{alg_rbanova} with $\Nmd$ samples drawn from the
prior distribution $\pi(\xi)$, and start the MCMC iterations with this initial model. 
Initializing a Markov chain $\Xi^*:=\{\xi^{(1)}\}$ where $\xi^{(1)}$ is a sample from the prior $\pi(\xi)$, 
for each $j\geq1$, we first draw a candidate sample $\xi^*$ from a proposal distribution which is denoted by $\pi(\cdot|\xi^{(j)})$,
and evaluate the system output corresponding to $\xi^{(j)}$ using Algorithm~\ref{alg_rbpre}, which is denoted by $G^r_{\pP}(\xi^{(j)})$. 
After that, a Metropolis acceptance ratio is computed through
\begin{eqnarray}
a:=\min\left(1, \frac{\pi_\epsilon\left(d-G^r_{\pP}\left(\xi^*\right)\right)\pi(\xi^*)}{\pi_\epsilon\left(d-G^r_{\pP}\left(\xi^{(j)}\right)\right)\pi\left(\xi^{(j)}\right)}
\frac{\pi\left(\xi^{(j)}|\xi^*\right)}{\pi\left(\xi^{*}|\xi^{(j)}\right)}\right).
\end{eqnarray}
With probability $a$, the candidate sample is accepted, i.e., $\xi^{(j)}:=\xi^*$; otherwise, the candidate sample is rejected, 
 i.e., $\xi^{(j)}:=\xi^{(j-1)}$. The Markov chain is then augmented with $\xi^{(j)}$, i.e., $\Xi^*=\Xi^*\cup \xi^{(j)}$.
 After $\Nmd$ posterior samples are generated, 
 the RB-ANOVA model is updated---the RB-ANOVA model data $\{c,\,\uh(x,c),\, \pP,\, \{Q_t\}_{t\in\pP}\}$ are reconstructed 
  using Algorithm \ref{alg_rbanova} with these $\Nmd$ posterior samples.
The MCMC procedure continues with the new RB-ANOVA model.
The RB-ANOVA model is reconstructed periodically every $\Nmd$ MCMC iterations, 
until certain stoping conditions are satisfied. 
Namely, as the number of MCMC samples increases, it is expected that the resulting RB-ANOVA model may not vary much. 
Thus, we terminate the reconstruction procedure  if the new model data 
and current model data are similar. 
Specifically, the index set $\pP$ is used to serve as the stoping criterion: 
the model reconstruction procedure is stopped  if the new and the current index sets are the same. 

This new adaptive RB-ANOVA-MCMC procedure is formally presented in Algorithm \ref{alg_rbanovamcmc}.
In the inputs of this algorithm, $N$ refers to the desired number of posterior samples to generate,
and $\Nmd$ is the sample size to generate the RB-ANOVA model. The variable {\em Update\_Label}
is used to label whether to stop updating the RB-ANOVA model during the MCMC iterations. 

Finally, we provide some discussions on how the use of the posterior distribution may improve the performance
of the model reduction. 
The improvement is two-fold: 
it improves the efficiency of both the ANOVA model (for the random parameters) and the reduced basis model (for the physical parameters). 
First, for the ANOVA model, both the anchor point and the important terms 
are selected based on some statistical moments of the random parameters. 
In particular, it has been discussed in \cite{gaohes10} that the efficiency of an ANOVA expansion depends critically on the 
choice of the anchor point---to achieve a given level of accuracy, a properly chosen anchor point can lead to 
a small number of expansion terms in \eqref{anova} or \eqref{rb_anova}, 
and they have suggested that an effective choice of the anchor point is the mean of the random parameters \cite{gaohes10}. 
Moreover, the active terms of the ANOVA model are also selected using the relative means. 
In a Bayesian problem, the random parameters are essentially distributed according to the posterior rather than the prior, and 
thus estimating these moments with respect to the posterior distribution should yield a much more accurate ANOVA representation 
than that with the prior. 
On the other hand, constructing the input sample set to generate the RB-ANOVA model  
from the posterior can also improve the performance of  the reduced basis model and  
the argument here is similar as that in \cite{li2014adaptive,cui2015data}: 
since the RB functions are chosen with respect to the input samples, 
constructing input samples from the posterior can ensure that the basis functions are mostly 
distributed in the high probability regions of the posterior, and the resulting RB model may be of higher accuracy in those regions. 
We will demonstrate that the proposed method can significantly improve the performance in  Section~\ref{section_numerical}.

\begin{algorithm}
\caption{The adaptive RB-ANOVA-MCMC algorithm}
\label{alg_rbanovamcmc}
\begin{algorithmic}[1]
\STATE Input: $N$ and $\Nmd$. 
\STATE Compute the RB-ANOVA model data $\left\{c, u^h(x,c), \pP, \{Q_t\}_{t\in\pP}\right\}$ 
using Algorithm \ref{alg_rbanova} with $\Nmd$ samples drawn from the prior distribution $\pi(\xi)$.
\STATE Draw a sample $\xi^{(1)}$ from the prior, and initialize the Markov chain $\Xi^* := \left\{\xi^{(1)}\right\}$. 
\STATE Let {\em Update\_Label} $:=1$.
\FOR{$j=1,\ldots,N-1$}
\STATE Draw $\-\xi^*\sim \pi\left(\cdot |\-\xi^{(j)}\right)$.
\STATE Compute the RB-ANOVA output $G^r_{\pP}\left(\xi^*\right)$ using Algorithm \ref{alg_rbpre}.
\STATE Compute the acceptance ratio 
\begin{equation*}
a = 
\min\left(1, \frac{\pi_\epsilon\left(d-G^r_{\pP}\left(\xi^*\right)\right)\pi(\xi^*)}{\pi_\epsilon\left(d-G^r_{\pP}\left(\xi^{(j)}\right)\right)\pi(\xi^{(j)})}
\frac{\pi\left(\xi^{(j)}|\xi^*\right)}{\pi\left(\xi^{*}|\xi^{(j)}\right)}\right).
\end{equation*}
\STATE Draw $\rho\sim U[0,1]$.
                  \IF{$ \rho<a$} 
  \STATE                  Let $\-\xi^{(j+1)} = \-\xi^*$, 
                  \ELSE 
                  \STATE Let $\-\xi^{(j+1)}=\-\xi^{(j)}$. 
                  \ENDIF
\IF{$j\mod \Nmd=0$  and {\em Update\_Label} $=1$}
\STATE Store the current ANOVA index set $\pP'=\pP$.
\STATE Update the RB-ANOVA model data $\left\{c, u^h(x,c), \pP, \{Q_t\}_{t\in\pP}\right\}$ 
using Algorithm \ref{alg_rbanova} with the last $\Nmd$ samples in the chain $\left\{\xi^{(j-\Nmd+1)},\ldots,\xi^{(j)}\right\}\subset \Xi^*$.
\IF{$\pP$ is the same as $ \pP'$}
\STATE Stop future RB-ANOVA model updates through setting {\em Update\_Label} $:=0$.
\ENDIF 
\ENDIF
\ENDFOR
\end{algorithmic}
\end{algorithm}

\section{Numerical study}\label{section_numerical}
\begin{figure}[!htp]
\centerline{
\begin{tabular}{cc}
\includegraphics[width=5.2cm,height=4.5cm]{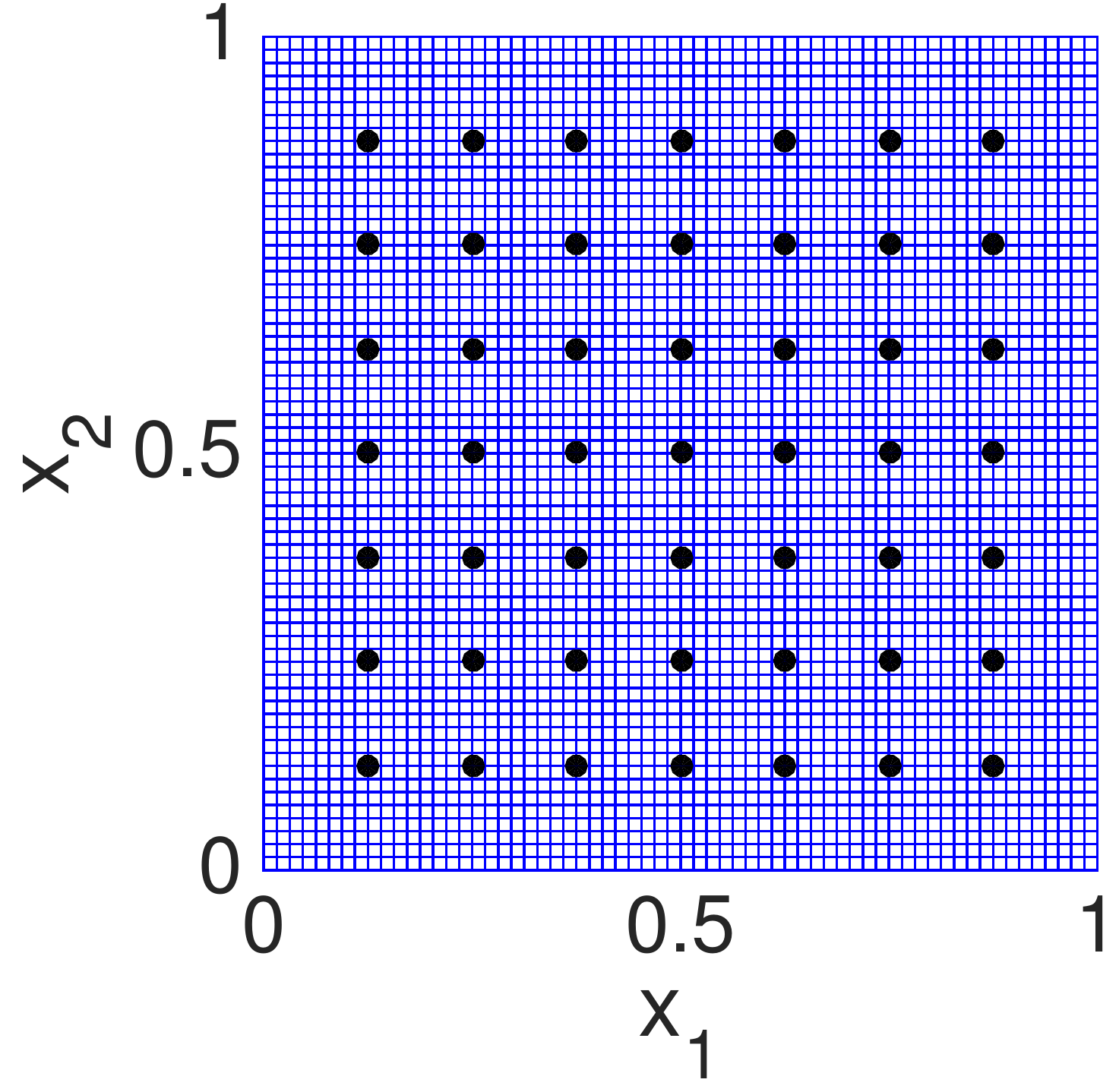}&
\includegraphics[width=5.5cm,height=4.5cm]{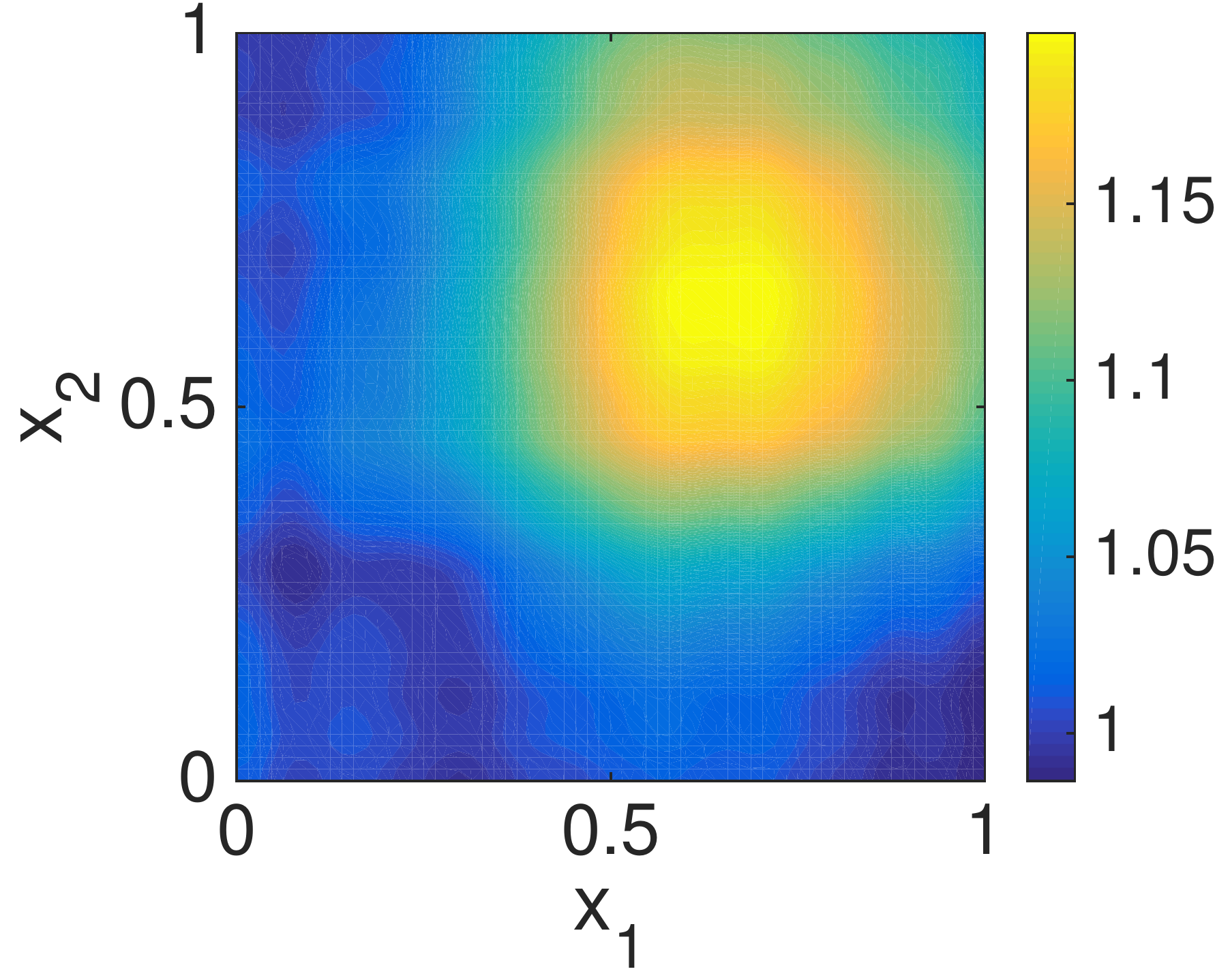}\\
(a)  FEM grids and sensors & (b) The actual permeability field \\ 
\end{tabular}}
\caption{Setup of the numerical test.}
\label{f:data}
\end{figure}

The numerical examples considered are steady flows in porous media. 
Letting $a(x,\xi)$ denote a unknown permeability field and $u(x,\xi)$ the pressure head,
we consider the following  diffusion equation, 
\begin{subequations} \label{e:diff}
\begin{align}
-\div \left(a\left(x,\xi\right)\nabla u\left(x,\xi\right)\right)=1 
& \quad \textrm{in} \quad  D\times I^M, \label{diff1} \\
u\left(x,\xi\right)=0                                  
& \quad\textrm{on} \quad \partial D\times I^M,\label{diff2}
\end{align}
\end{subequations}
where $D\subset \dsR^{2}$ and the dimension of the parameter $M$ is specified when we parameterize 
the permeability field next. 
Given a realization of $\xi$, 
defining $H^1(D):=\{u: D \to \dsR, \int_{D} u^2\, {\rm d} D<\infty, \int_{D} (\partial u/ \partial x_1)^2\, {\rm d} D<\infty \textrm{ and }  \int_{D} (\partial u/ \partial x_2)^2\, {\rm d} D<\infty\}$
and $H_0^1(D):=\{v\in H^1(D)\,| \,  v=0 \textrm{ on } \partial D_D\}$, 
the weak form of \eqref{e:diff} is to find
$u(x,\xi)\in H_0^1(D)$ such that
$(a \nabla u, \nabla v) = (1,v)$ for all $v\in H_0^1(D)$.
We discretize in space using a bilinear finite element approximation \cite{brae97,elman05}.
The  spatial domain in the following numerical studies is taken to be $D=(0,1)\times(0,1)$. 
The problem is discretized in space on a uniform $65\times 65$ grid (the number of the spatial degrees 
of freedom is $N_h=4225$).
Our deterministic forward model $G(\xi)$
is defined to be a set collecting solution values
corresponding to measurement sensors---$\{u(x,\xi), \xi \in \md\}$ where the sensor set $\md$ 
in this work is defined to be the tensor product $\{x_i\}\otimes\{y_i\}$ of the  one-dimensional grids:
$x_i=0.125i$, $y_i=0.125i$, for $ i=1,\ldots, 7.$ 
We set the measurement noise $\bolde$ in \eqref{eq:setting} to
independent and identically distributed Gaussian distributions with mean  zero and standard deviation $0.001$.
Figure \ref{f:data} shows locations of sensors with the finite element grids
and the true permeability field used to generate the test data.

We parameterize the permeability field $a(x,\xi)$ 
by a truncated Karhunen--Lo\`eve (KL) expansion \cite{ghaspa03,babuska1,elmmil11} of a random field with 
mean function $a_0(x)$, standard deviation $\sigma$ and covariance function
\begin{eqnarray}
 Cov(x,y)=\sigma^2 \exp\left(-\frac{|x_1-y_1|}{\covl}-\frac{|x_2-y_2|}{\covl}\right), \label{covariance}
\end{eqnarray}
where $\covl$ is the correlation length.
The KL expansion is expressed as
\begin{eqnarray}
a(x,\xi)=a_0(x)+\sum_{k=1}^M\sqrt{\lambda_k}a_k(x)\xi_k, \label{kl}
\end{eqnarray}
where $\{\xi_k\}^M_{k=1}$ are random variables,
$M$ is the number of KL modes retained, 
$a_k(x)$ and $\lambda_k$ are the eigenfunctions and eigenvalues of \eqref{covariance}.
We set $a_0(x)=1$ and $\sigma=0.25$ in the numerical studies. 
The priori distributions of $\{\xi_k\}^M_{k=1}$ are set to  
independent uniform distributions with range $I=[-1,1]$. 
Different values of the correlation length $\covl$ are studied.
As usual, we set $M$ large enough, such that $95\%$ of the total variance of the exponential covariance function
are captured \cite{powelm09}.

\subsection{The impact of priors}
Different priors are tested for this problem and we shall see how the priors affect the inference results. 
We specifically test the prior permeability fields  associated with 
four different values of the correlation length $\covl$ in \eqref{covariance}: $5$, $5/2$, $5/4$ and $5/8$. 
To capture  $95\%$ of the total variance of the covariance function, we set the number of KL modes retained (the dimension the 
patermeter $\xi$) as:
$M=4$ for $\covl=5$, $M=8$ for $\covl=5/2$, $M=23$ for $\covl=5/4$ and $M=73$ for $\covl=5/8$.

To generate posterior samples for comparison,
the MCMC method described in Algorithm~\ref{alg:mh} is first performed with the forward model  
evaluated by the finite element method, which is referred to as the full MCMC method.
We here draw $N=10^6$ posterior samples using full MCMC with each of the above four priors.
In all our numerical tests, the proposal distribution $\pi(\cdot|\xi^{(j)})$  
on line $3$ of Algorithm~\ref{alg:mh} (and on line $6$ of Algorithm~\ref{alg_rbanovamcmc}) 
is set to a multivariate Gaussian distribution with mean $\xi^{(j)}$ and covariance matrix $0.03^2I$, where $\xi^{(j)}$ is the $j$-th sample
in the Markov chain and $I\in\dsR^{M\times M}$ is an identity matrix. 
The acceptance rates (numbers of accepted samples divided by the total sample size)  
are $47\%$,$46\%$,$42\%$ and $26\%$ for $M=4$, $M=8$, $M=23$ and $M=73$ respectively, 
which indicates that the proposal is properly chosen \cite{roror01}. In addition as expected, the acceptance rate decreases as
the parameter dimension increases.

\begin{figure}[!htp]
\centerline{
\begin{tabular}{cc}
\includegraphics[width=5.5cm,height=4.5cm]{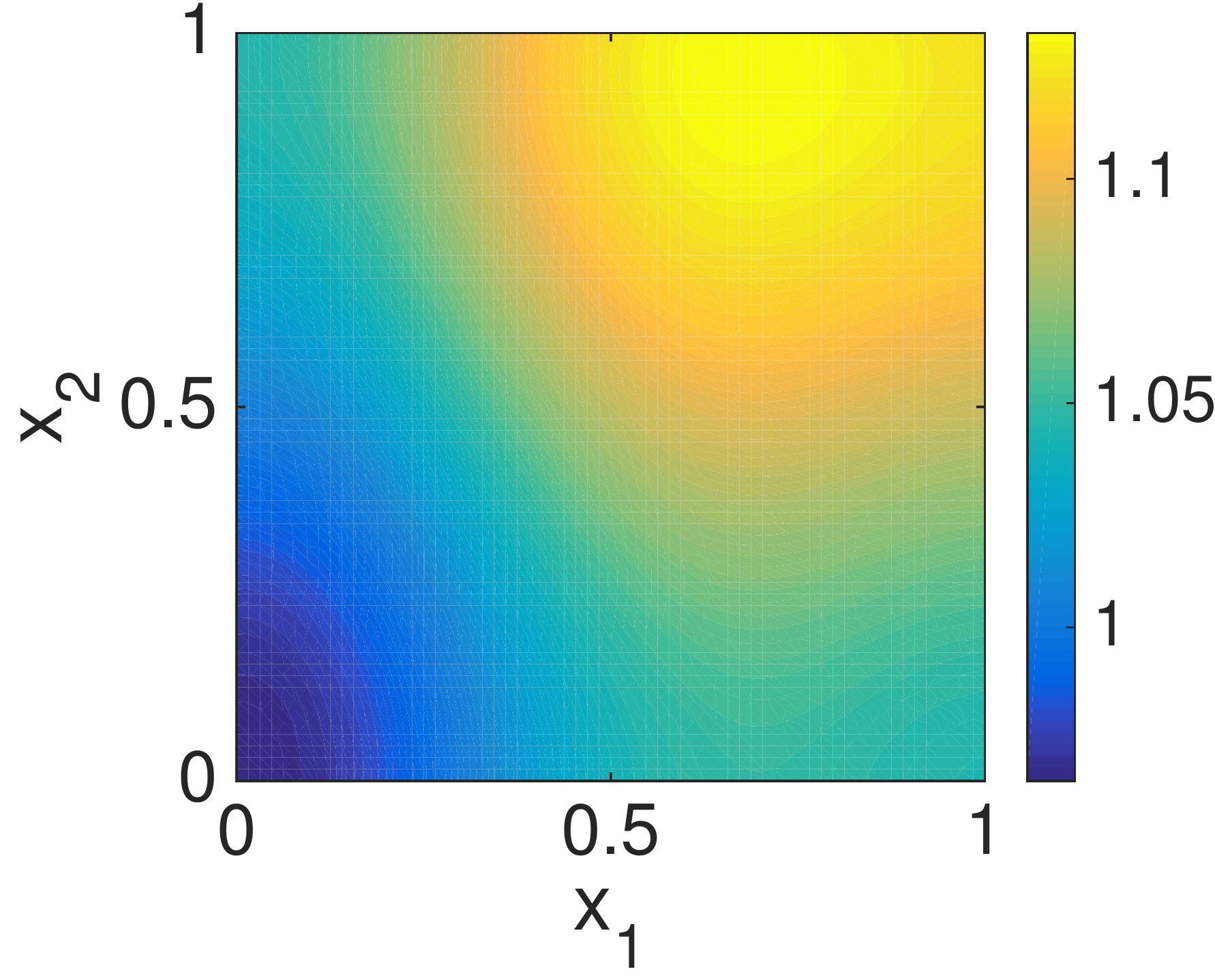} & 
\includegraphics[width=5.5cm,height=4.5cm]{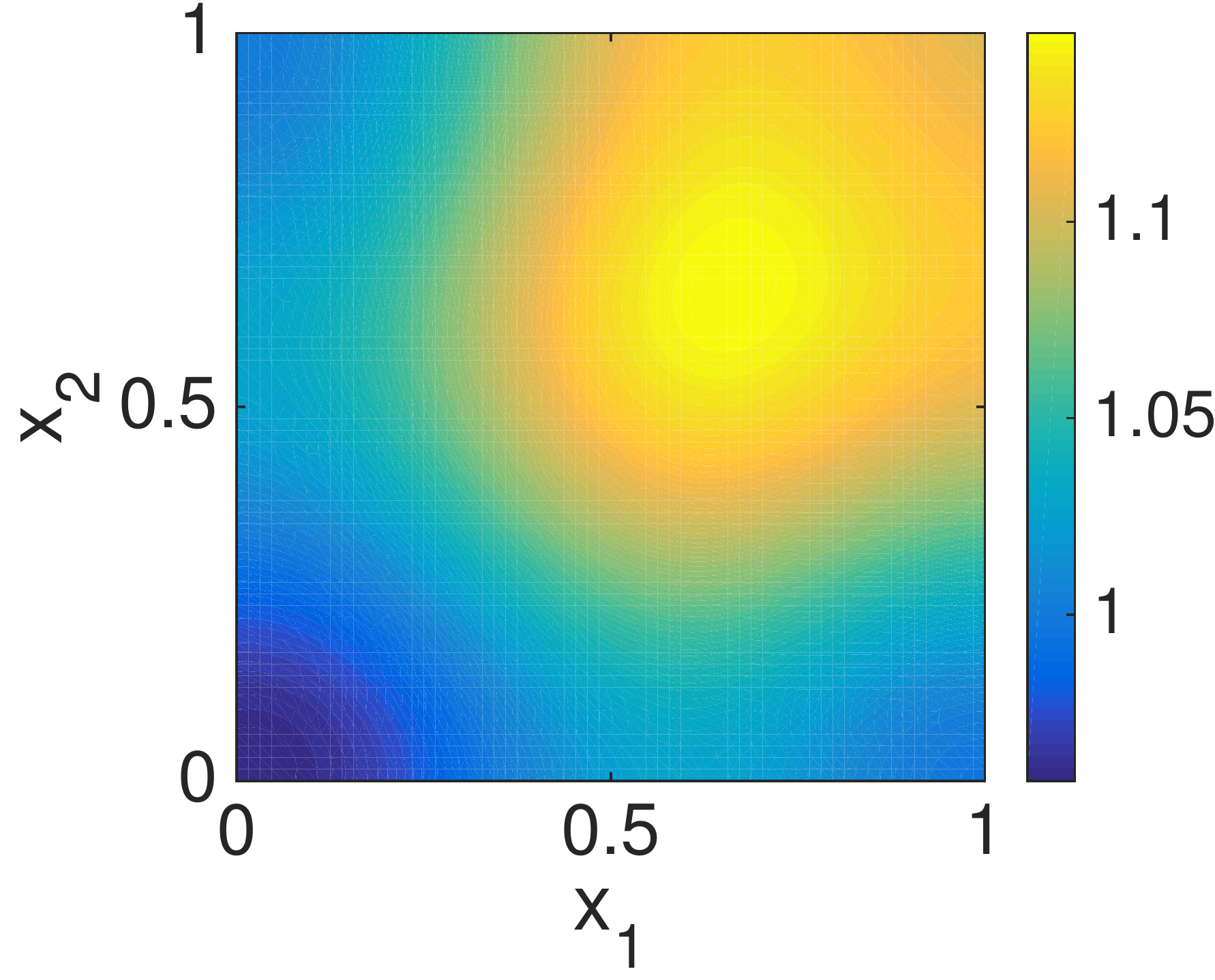} \\
(a) $L=5$, $M=4$ & (b) $L=5/2$, $M=8$ \\
\includegraphics[width=5.5cm,height=4.5cm]{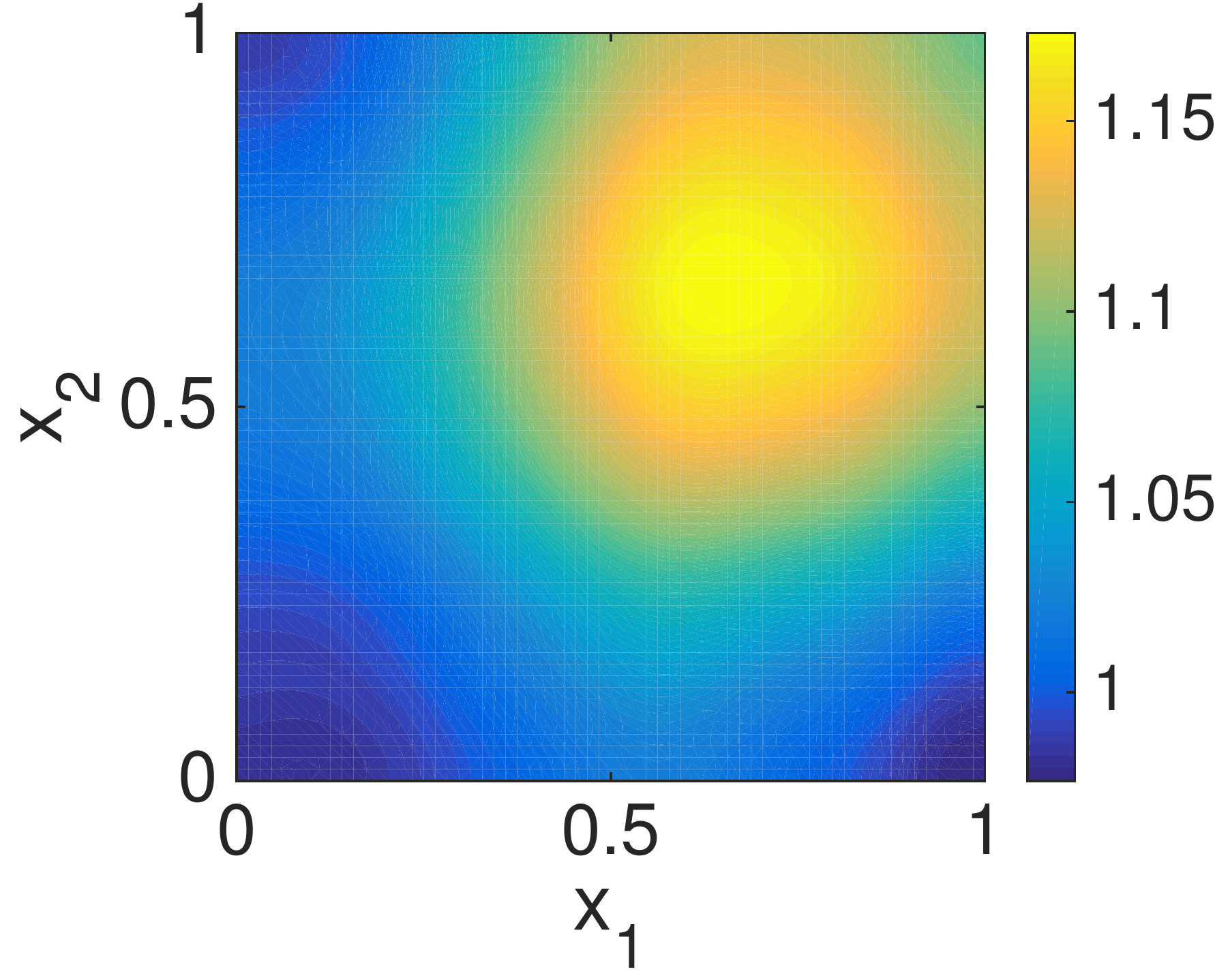} & 
\includegraphics[width=5.5cm,height=4.5cm]{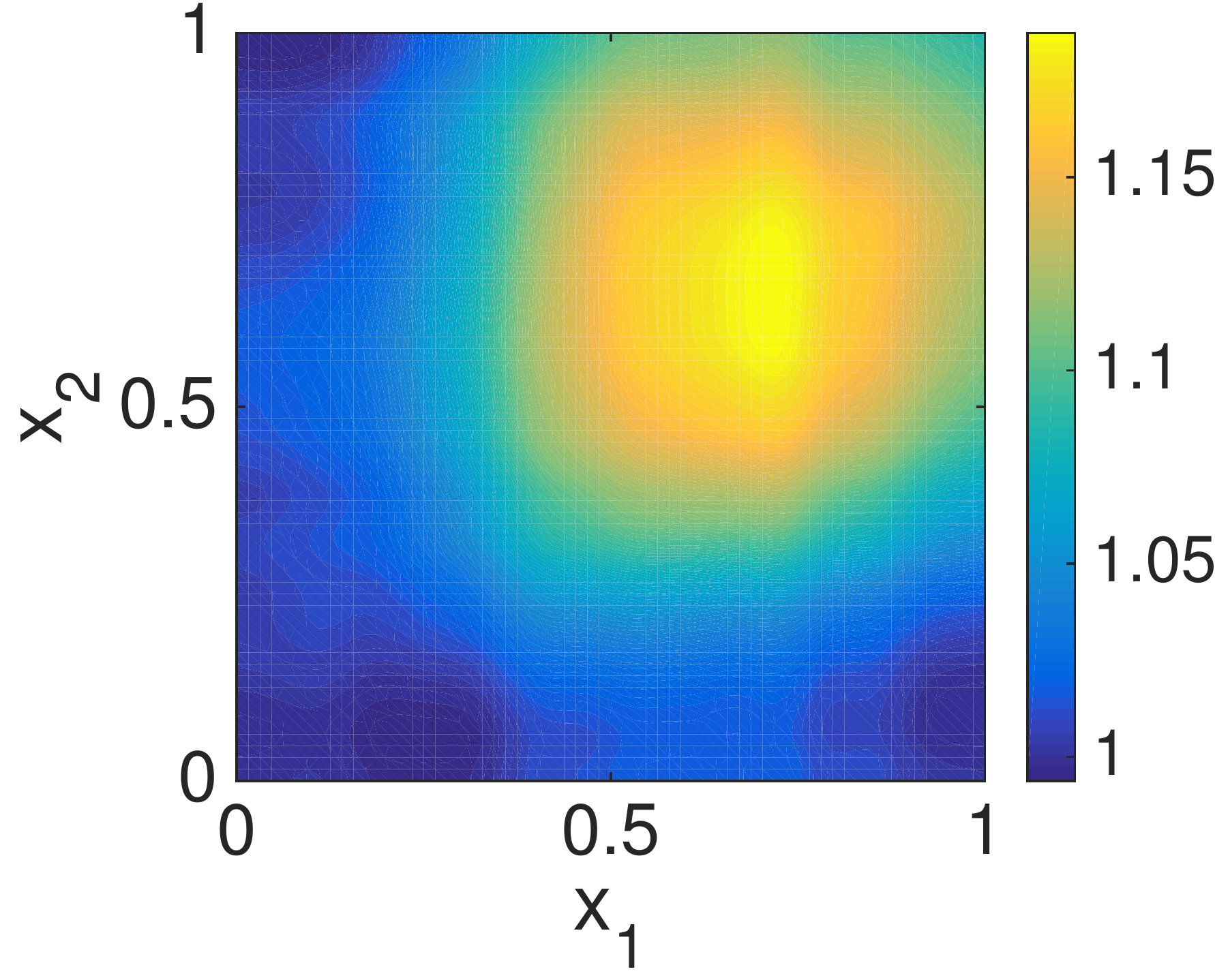} \\
(d) $L=5/4$, $M=23$ & (d) $L=5/8$, $M=73$ \\
\end{tabular}}
\caption{Full MCMC results.}
\label{f:full_mcmc}
\end{figure}

\begin{figure}[!htp]
\centerline{
\begin{tabular}{c}
\includegraphics[width=5cm,height=4cm]{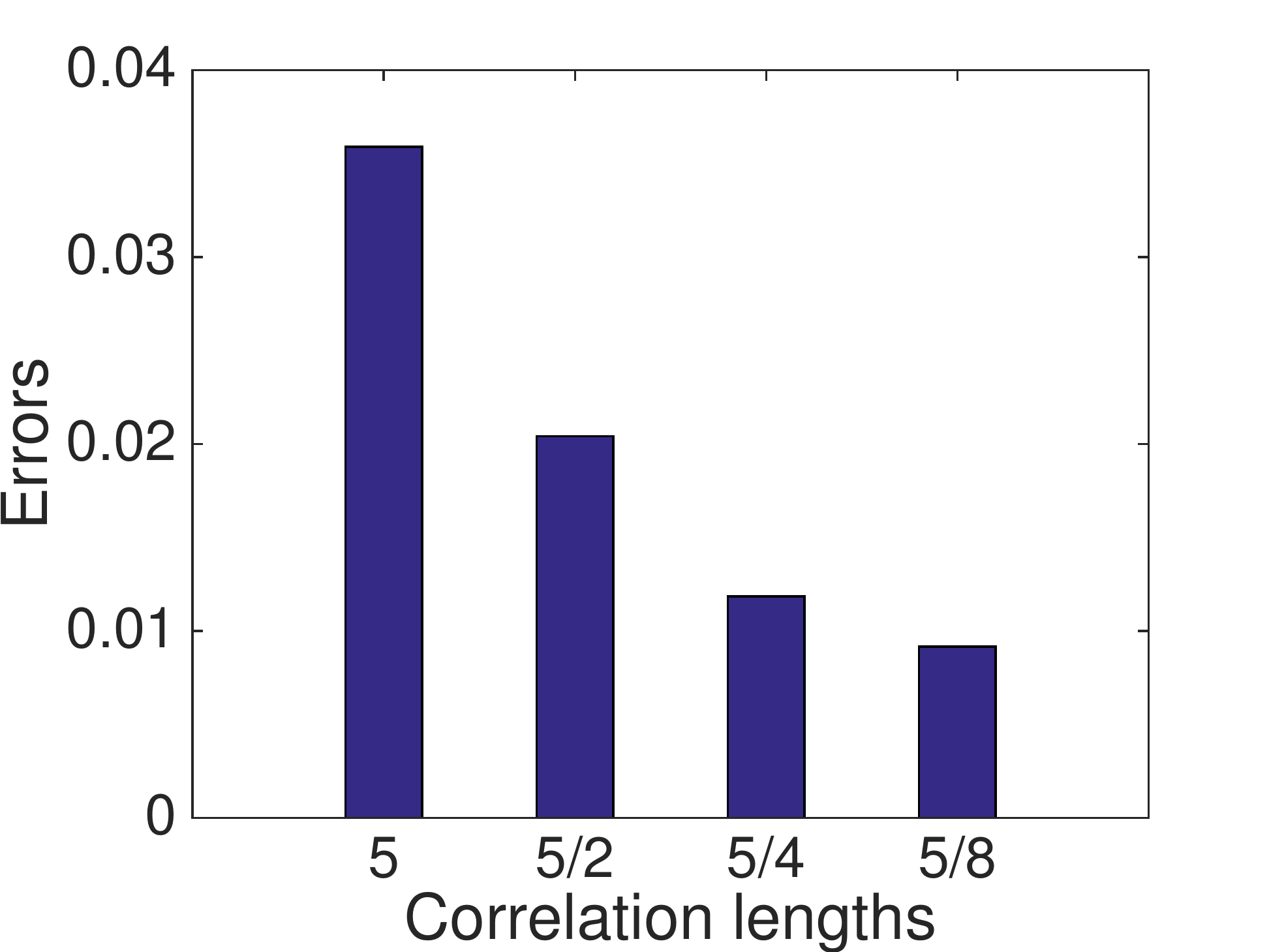} 
\end{tabular}}
\caption{Errors ($\epsilon_{\Xi^*}$) w.r.t. correlation lengths ($\covl$).}
\label{f:full_error}
\end{figure}

Figure \ref{f:full_mcmc} shows the estimated posterior mean permeability fields generated by the sample means of 
full MCMC, each of which is defined as 
\begin{eqnarray}
\pE_{\Xi^*} \left(a\left(x,\xi\right) \right)&:=&\sum_{\xi\in\Xi^*}\frac{a\left(x,\xi\right)}{|\Xi^*|},\label{mc_mean_field}
\end{eqnarray}
where $\Xi^*$ is the set of MCMC samples and $|\Xi^*|$ is its size.
It is clear that, as the correlation length $\covl$ reduces (the dimension of the parameter $M$ increases),
the estimated mean permeability field becomes visually similar to the actual  field shown in Figure~\ref{f:data}(b). 
In particular, for a large correlation length $\covl=5$, while the prior is very smooth, 
the estimated posterior mean permeability is also too smooth compared with the actual field. 
For a smaller correlation, e.g., $\covl=5/8$, the prior becomes less smooth, and 
the estimated  posterior mean permeability becomes more accurate. 
To assess the accuracy of the estimated posterior mean  permeability, 
we introduce the following quantity of errors
\begin{eqnarray}
\epsilon_{\Xi^*}:=\left\|\pE_{\Xi^*} \left(a\left(x,\xi\right) \right)-a_{\rm actual}\right\|_0 
\left/ 
\left\|a_{\rm actual}\right\|_0 \right., 
\label{e_mean_field}
\end{eqnarray}
where $a_{\rm actual}$ is the actual permeability field shown in Figure~\ref{f:data}(b).
Figure~\ref{f:full_error} shows the errors with respect to the correlation lengths, where it is clear that 
small correlation lengths lead to small errors for our test problem. This motivates us to focus on 
priors with small correlation lengths, which require high-dimensional parameterization.

\subsection{Performance of RB-ANOVA surrogates}
\begin{figure}[!htp]
\centerline{
\begin{tabular}{cc}
\includegraphics[width=5.5cm,height=4.5cm]{Fig/full_field_C125-eps-converted-to.pdf} & 
\includegraphics[width=5.5cm,height=4.5cm]{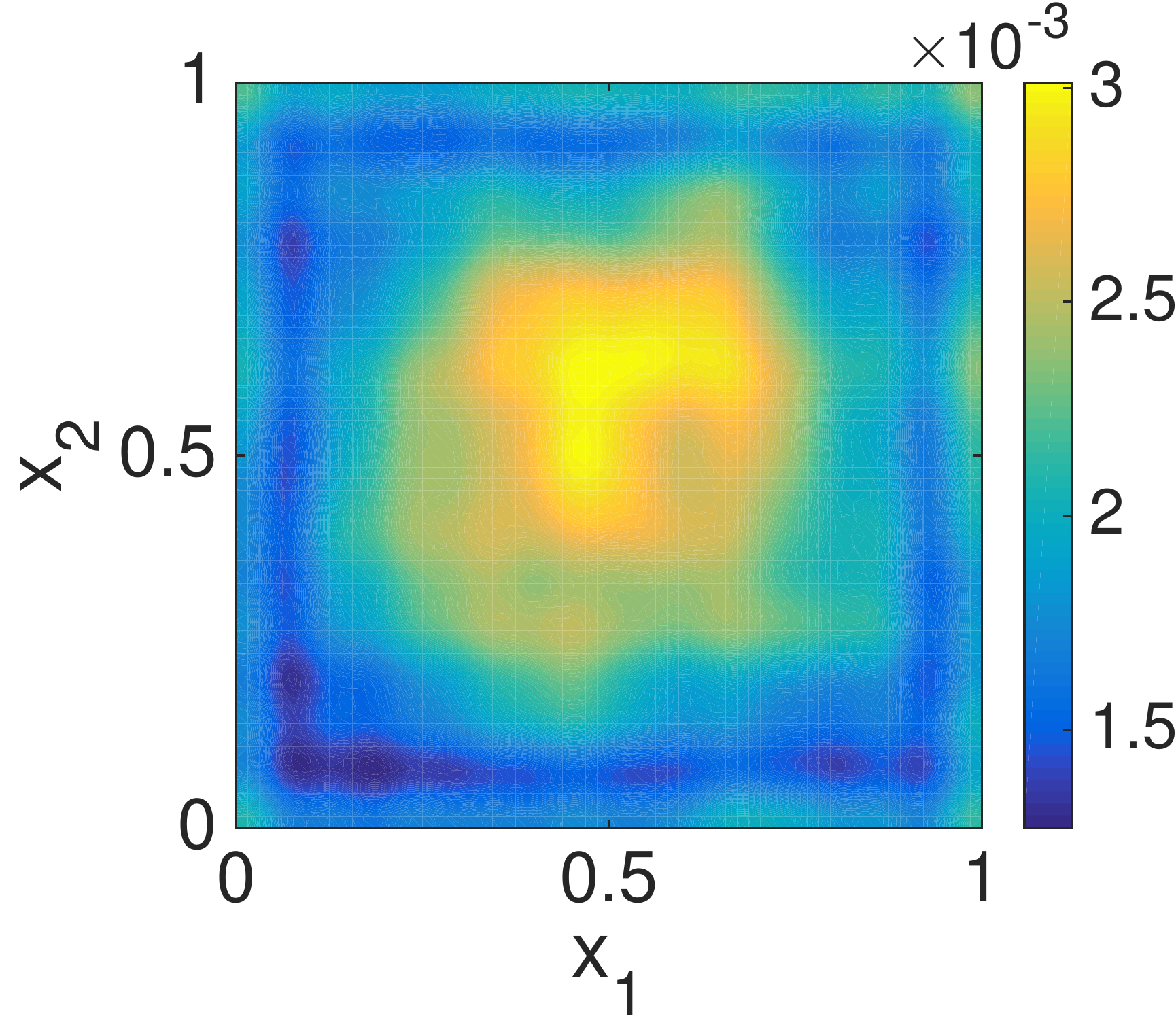} \\
(a) Mean, full MCMC & (b) Variance, full MCMC \\
\includegraphics[width=5.5cm,height=4.5cm]{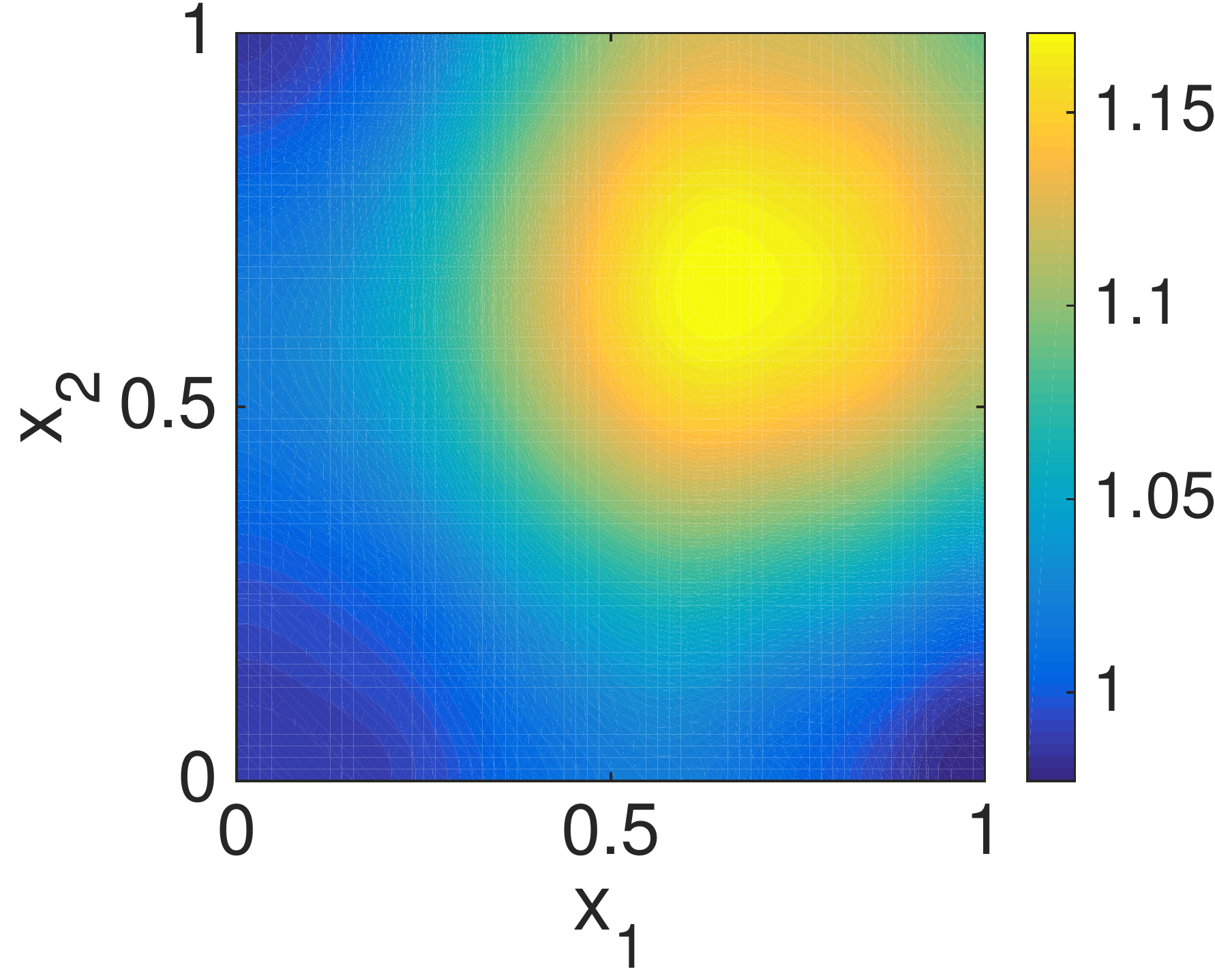} & 
\includegraphics[width=5.5cm,height=4.5cm]{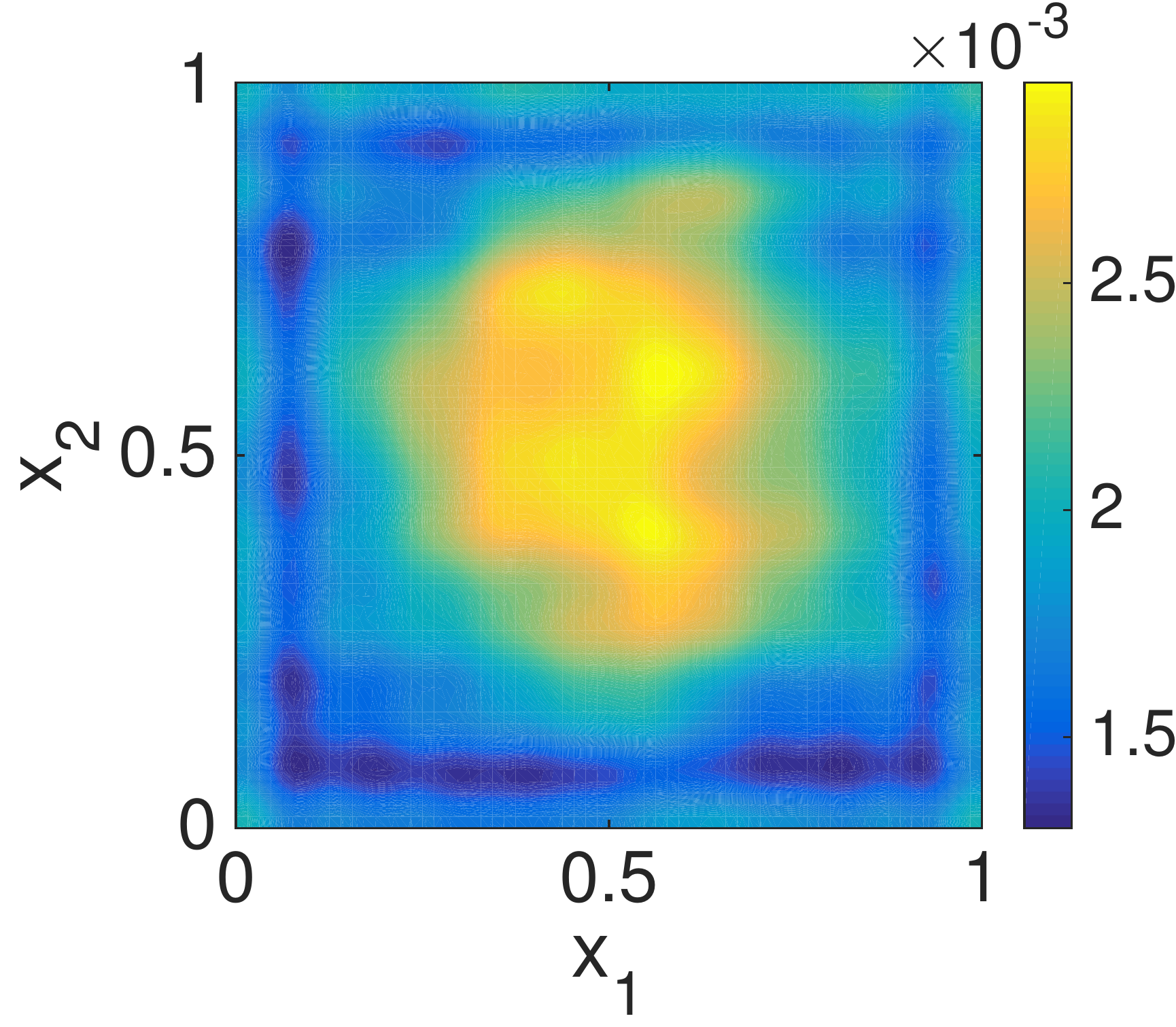} \\
(c) Mean, prior RB-ANOVA-MCMC & (d) Variance, prior RB-ANOVA-MCMC\\
\includegraphics[width=5.5cm,height=4.5cm]{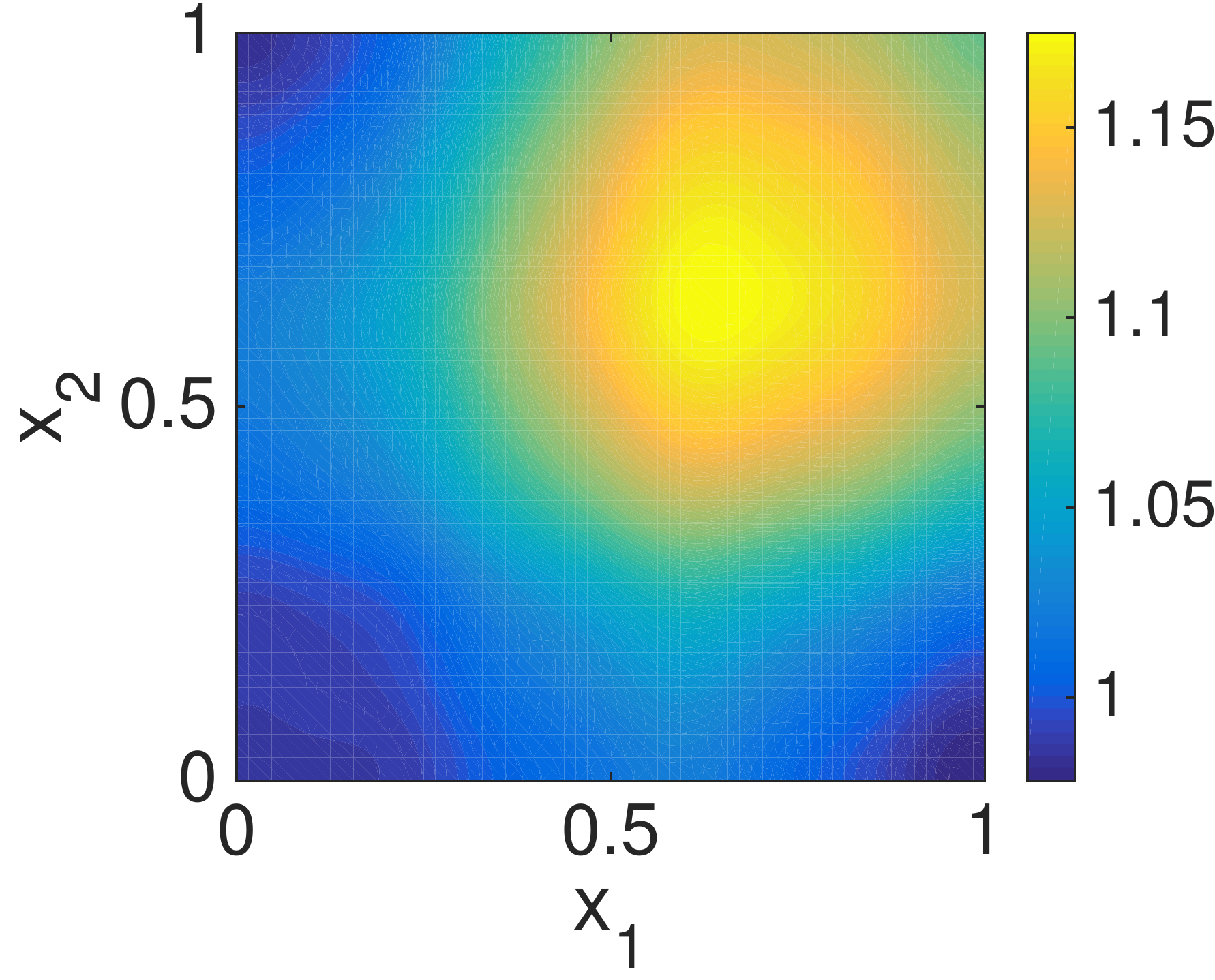} & 
\includegraphics[width=5.5cm,height=4.5cm]{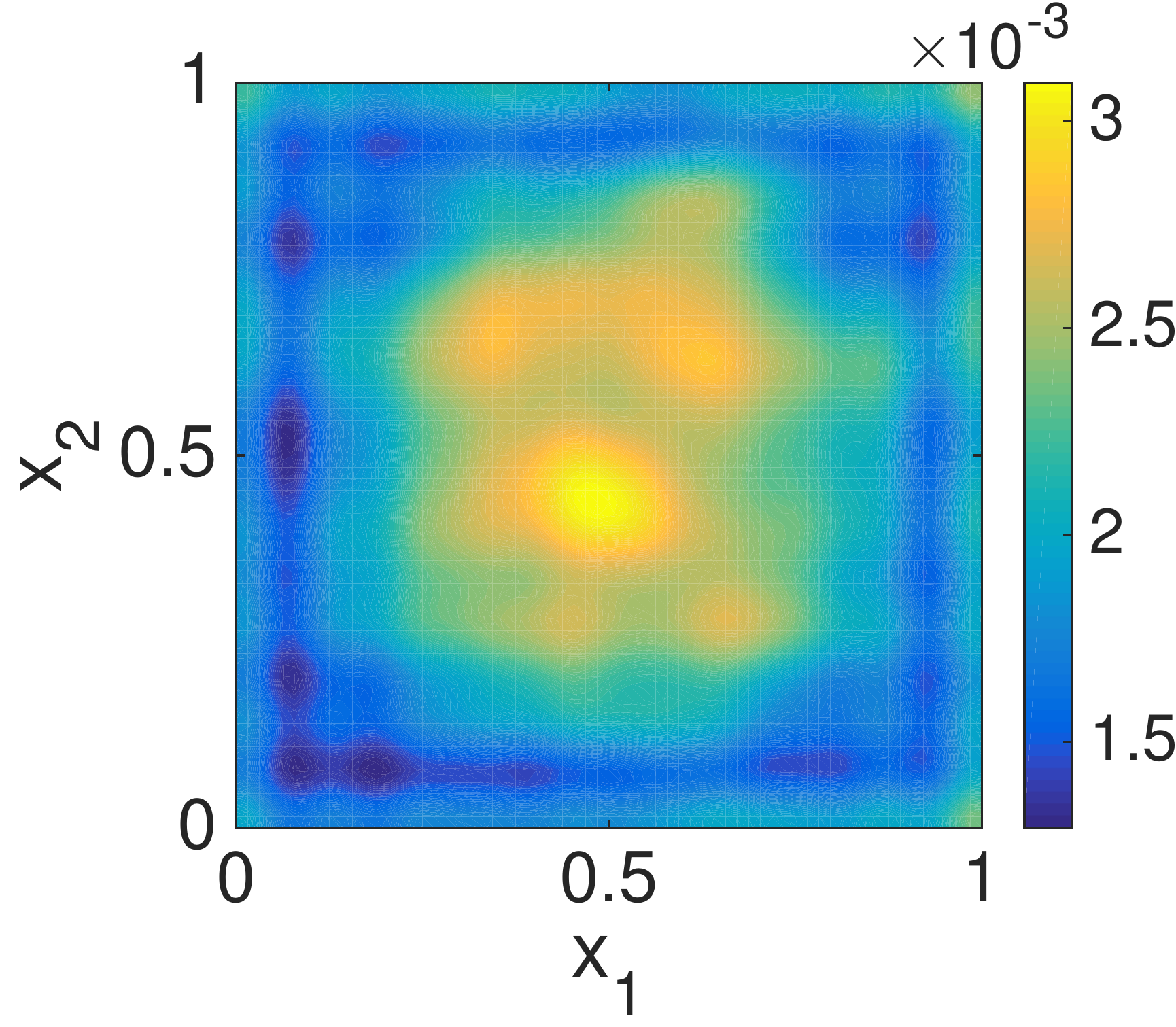} \\
(e) Mean, adaptive RB-ANOVA-MCMC & (f) Variance, adaptive RB-ANOVA-MCMC\\
\end{tabular}}
\caption{Estimated mean and variance fields for $\covl=5/4$ with $M=23$.}
\label{f:rb_mcmc_field_C125}
\end{figure}

We here focus on the two high-dimensional cases in our test problem ($\covl=5/4$ with 
$M=23$ and $\covl=5/8$ with  $M=73$), and test the RB-ANOVA-MCMC approach for these two cases.
For comparison, an unadaptive version of RB-ANOVA-MCMC is also tested in addition to the 
the adaptive RB-ANOVA-MCMC (Algorithm~\ref{alg_rbanovamcmc}). The unadaptive version,  which is referred to as
the prior RB-ANOVA-MCMC method in the following, uses samples from 
the prior distribution to generate the RB-ANOVA model through Algorithm~\ref{alg_rbanova} and performs the MCMC iterations
using this model. 
There are three tolerance parameters that need to be specified for generating the RB-ANOVA model in Algorithm~\ref{alg_rbanova}:
$tol_{pod}$ for selecting singular vectors in POD on line 9 (details are discussed in Section~\ref{section_rb}), 
$tol_{rb}$ on line 13 and $tol_{anova}$ on line 24. 
Following the discussion of our work \cite{liaolin16}, we set them all to $10^{-4}$ in this work.
For both prior and adaptive RB-ANOVA-MCMC, $10^3$ samples are used to generate the RB-ANOVA model,
i.e., $\Nmd=10^3$ in Algorithm~\ref{alg_rbanovamcmc}. 
Figure~\ref{f:rb_mcmc_field_C125} shows estimated mean and variance fields
for the case $\covl=5/4$ with $M=23$, 
generated by the three approaches: full MCMC, prior RB-ANOVA-MCMC, and adaptive RB-ANOVA-MCMC
respectively  with $10^6$ samples. Here, the estimated mean fields are computed through \eqref{mc_mean_field}, and 
the estimated variance fields are computed through
\begin{eqnarray}
\pV_{\Xi^*} \left(a\left(x,\xi\right) \right)&:=&\sum_{\xi\in\Xi^*}\frac{1}{|\Xi^*|}\bigg(a\left(x,\xi\right)-\pE_{\Xi^*} \left(a\left(x,\xi\right) \right)\bigg)^2 \label{mc_var_field},
\end{eqnarray}
where $\pE_{\Xi^*} (a(x,\xi) )$ is defined in \eqref{mc_mean_field} and $\Xi^*$ is the posterior sample set generated by each 
of the three approaches. From Figure~\ref{f:rb_mcmc_field_C125}, the estimated mean and variance fields 
generated by prior and adaptive RB-ANOVA-MCMC look very similar to those generated by full MCMC. For the case $\covl=5/8$
with $M=73$, Figure~\ref{f:rb_mcmc_field_C0625} shows that the estimated mean and variance fields generated by the three 
approaches are also very similar. 

\begin{figure}[!htp]
\centerline{
\begin{tabular}{cc}
\includegraphics[width=5.5cm,height=4.5cm]{Fig/full_field_C0625-eps-converted-to.pdf} & 
\includegraphics[width=5.5cm,height=4.5cm]{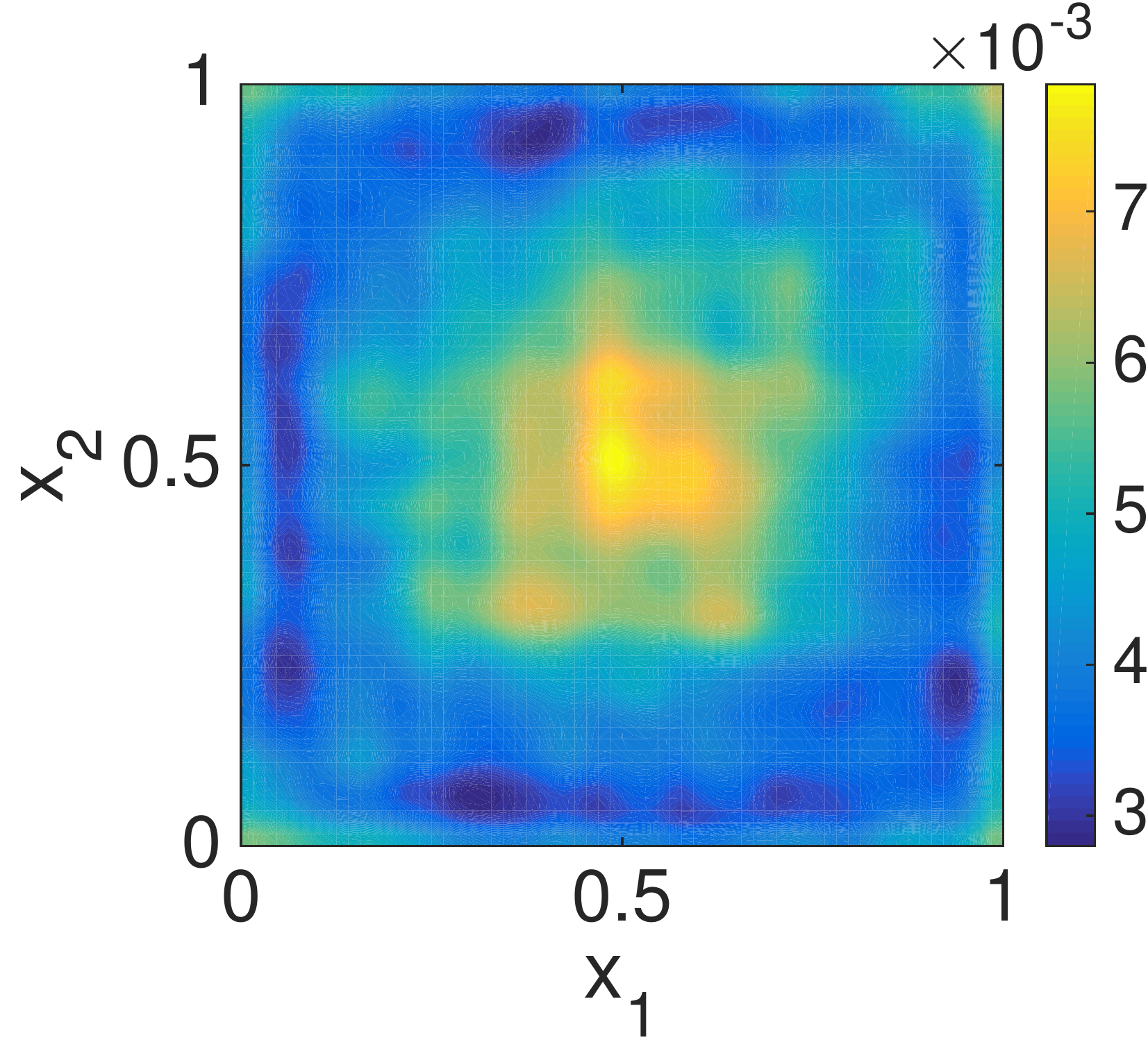} \\
(a) Mean, full MCMC & (b) Variance, full MCMC \\
\includegraphics[width=5.5cm,height=4.5cm]{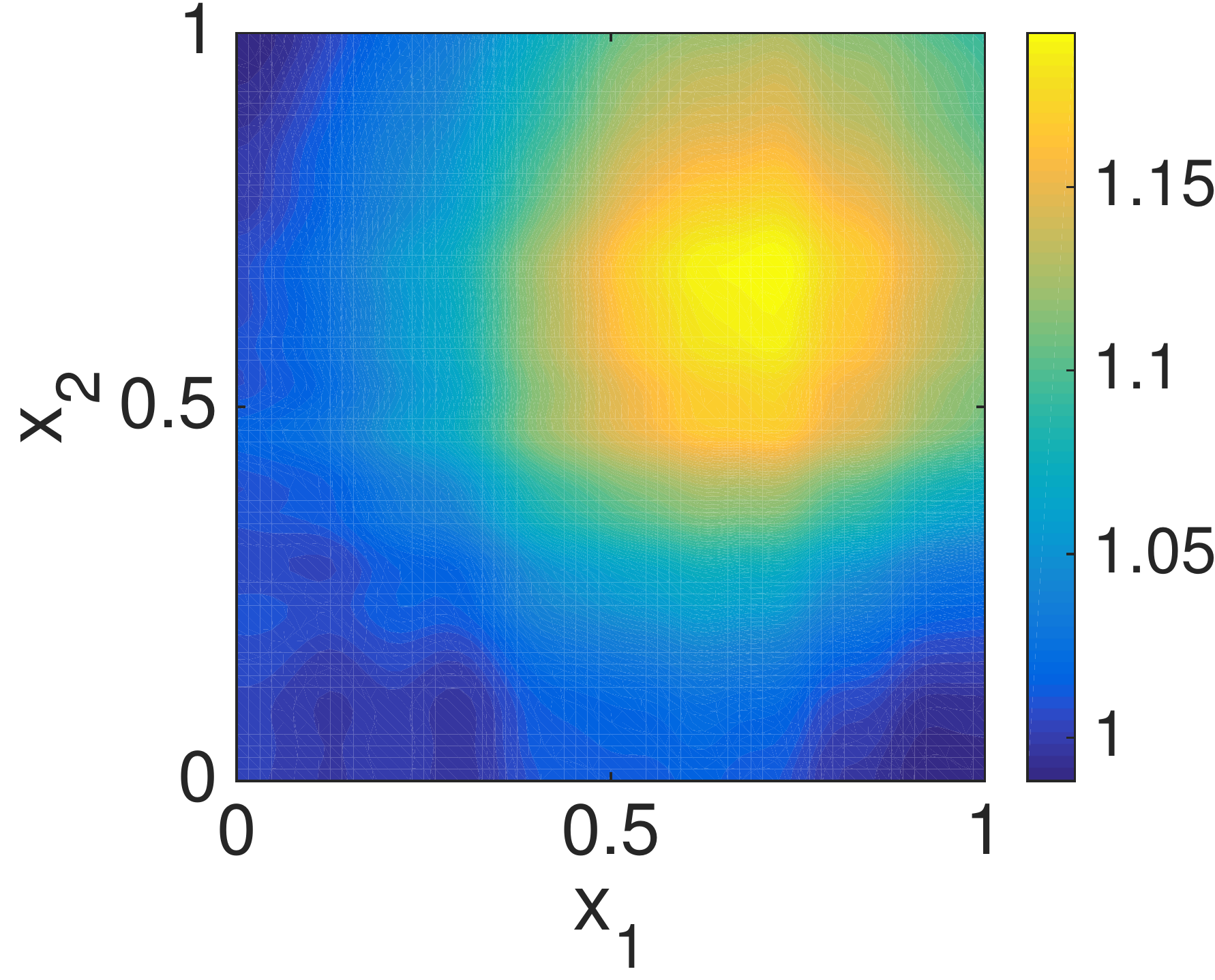} & 
\includegraphics[width=5.5cm,height=4.5cm]{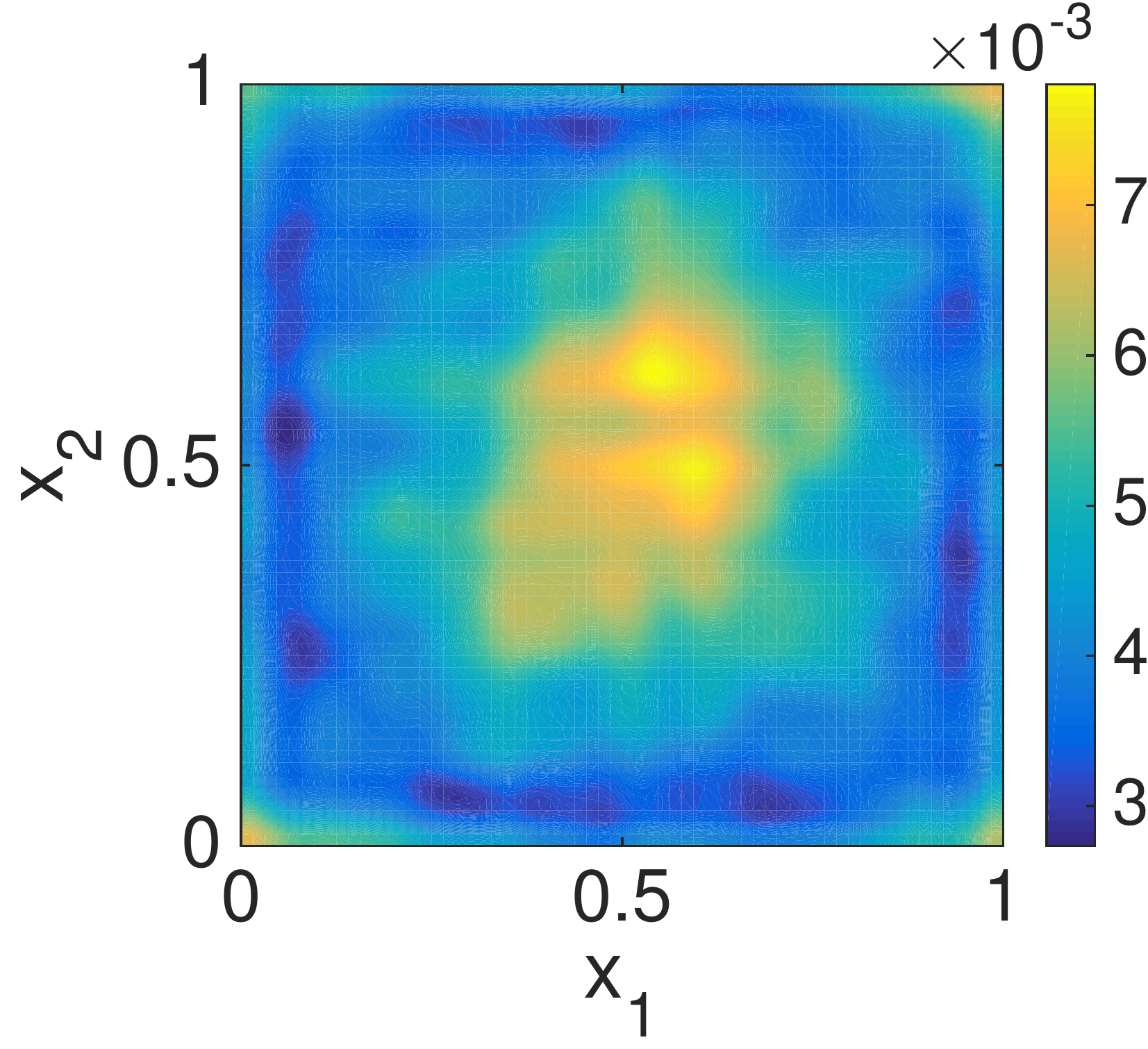} \\
(c) Mean, prior RB-ANOVA-MCMC & (d) Variance, prior RB-ANOVA-MCMC\\
\includegraphics[width=5.5cm,height=4.5cm]{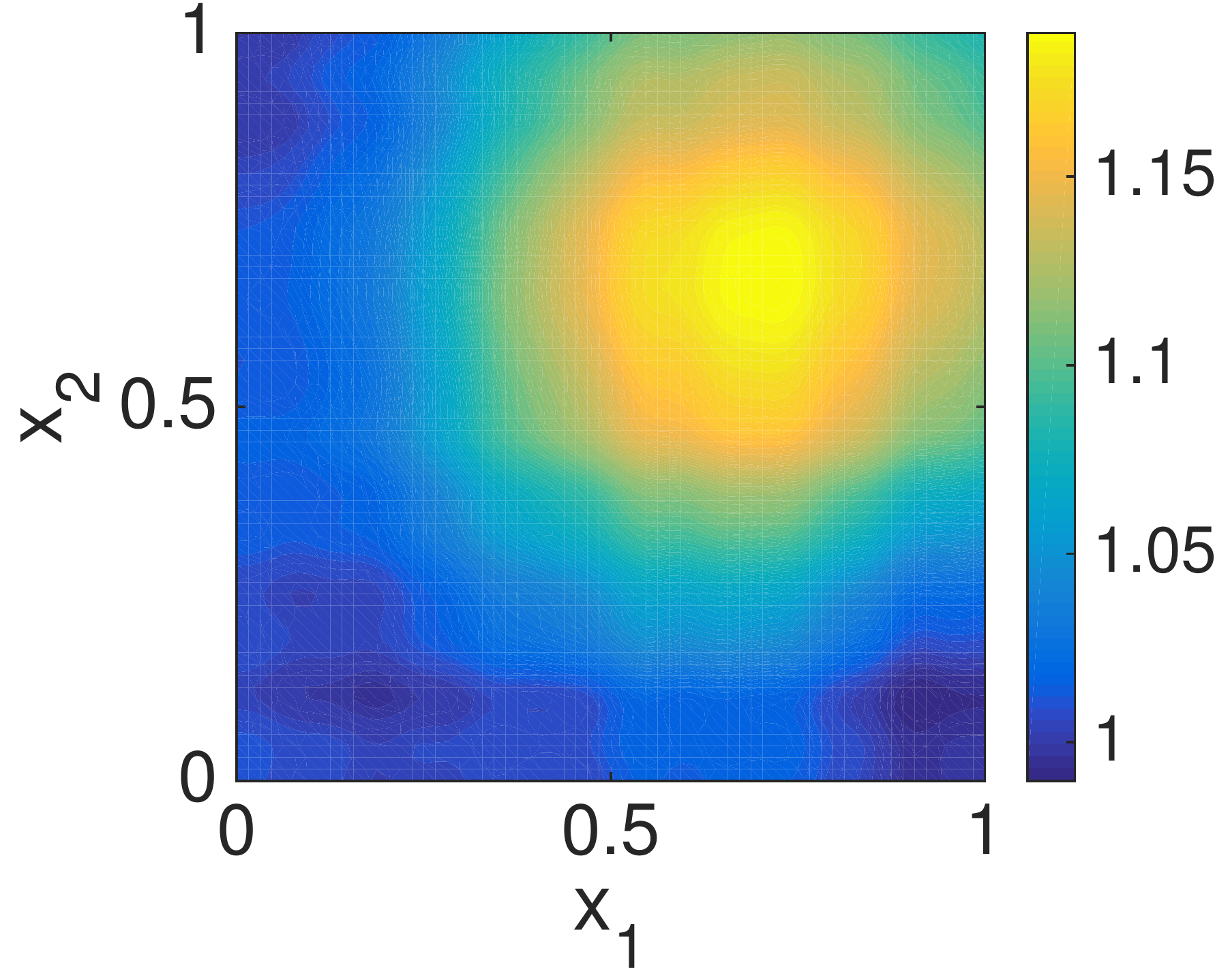} & 
\includegraphics[width=5.5cm,height=4.5cm]{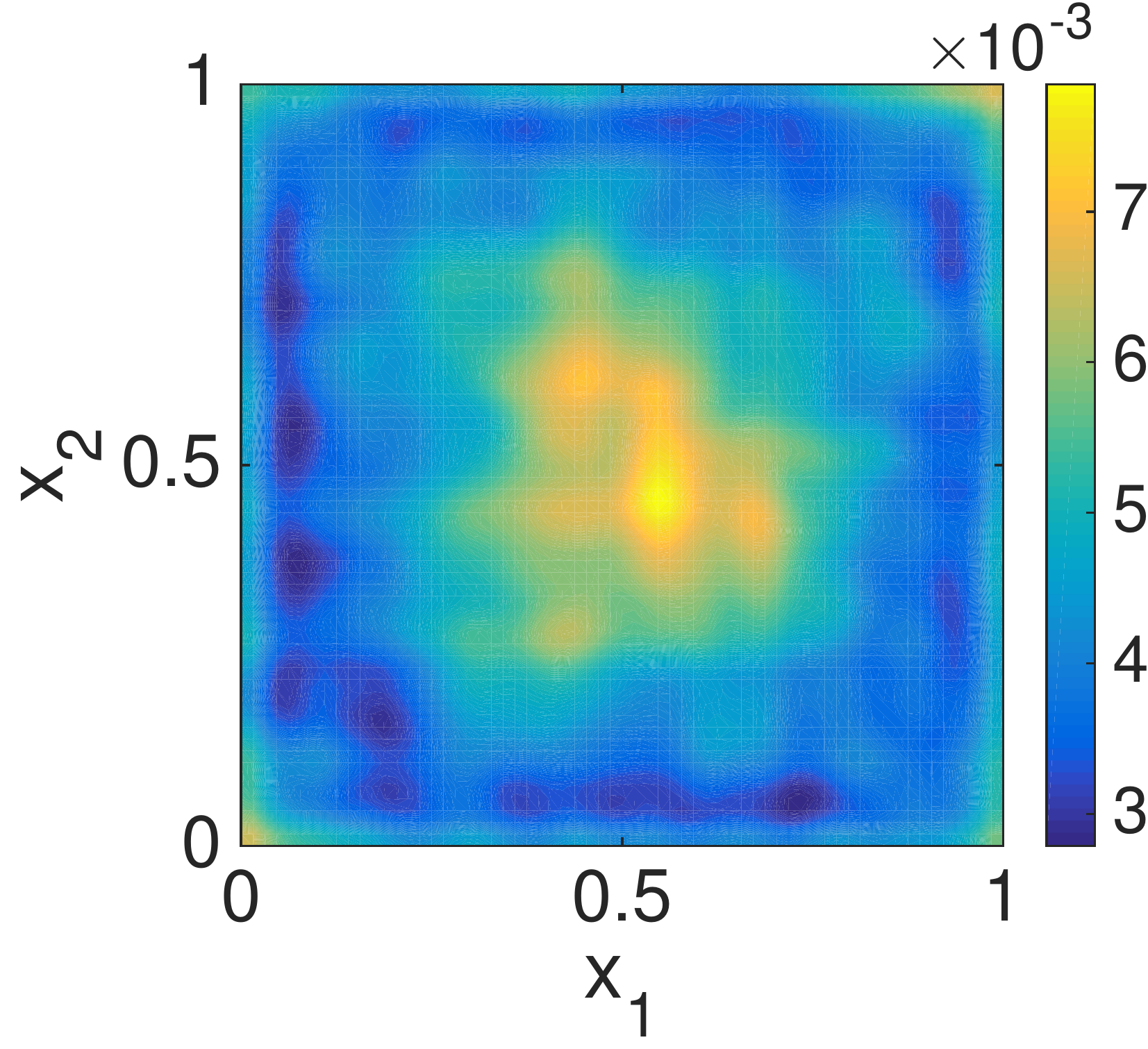} \\
(e) Mean, adaptive RB-ANOVA-MCMC & (f) Variance, adaptive RB-ANOVA-MCMC\\
\end{tabular}}
\caption{Estimated mean and variance fields for $\covl=5/8$ with $M=73$.}
\label{f:rb_mcmc_field_C0625-eps-converted-to.pdf}
\end{figure}

As discussed in Section~\ref{section_intro}, the main cost of the MCMC procedure comes 
from evaluating the forward model. For full MCMC, the forward model is evaluated using the finite element method,
while it is evaluated using the RB-ANOVA model in our RB-ANOVA-MCMC approach.
To assess the costs, we adopt the computational cost model for reduced basis methods developed in our recent work \cite{liaolin16},
which is based on counting relative sizes of linear systems (algebraic versions of (\ref{fem}) and (\ref{rb})). 
In this cost model, for a given finite element degrees of freedom $N_h$, the cost for solving a full system \eqref{fem} 
is defined to be a {\em cost unit}, which is assumed to be independent of the parameter $\xi$. 
The cost of solving a reduced problem (\ref{rb}) with size $N_r$ is modelled by $N_r/N_h$.
So, the cost of full MCMC is the number of forward model evaluations (see Algorithm~\ref{alg:mh}), 
and the cost of our adaptive RB-ANOVA-MCMC is the sum of the costs for solving reduced systems  (\ref{rb}) and 
full systems \eqref{fem} involved Algorithm~\ref{alg_rbanovamcmc}. In addition, it is clear that the cost of prior 
RB-ANOVA-MCMC is the sum of the costs in the construction procedure (Algorithm~\ref{alg_rbanova}) and the costs
of using  Algorithm~\ref{alg_rbpre} to evaluate forward models in the MCMC iterations.

\begin{figure}[!htp]
\centerline{
\begin{tabular}{cc}
\includegraphics[width=5.5cm,height=4.5cm]{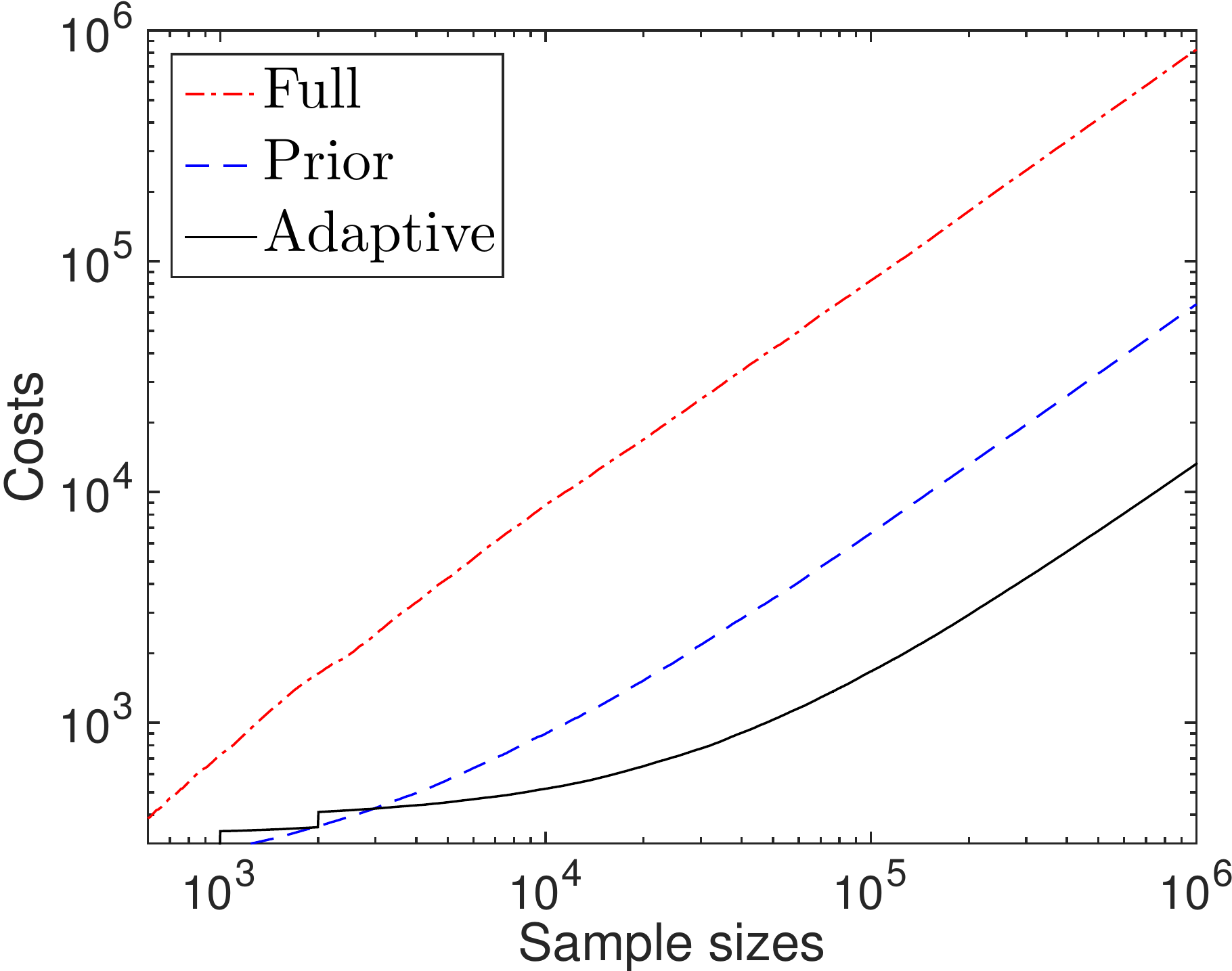} & 
\includegraphics[width=5.5cm,height=4.5cm]{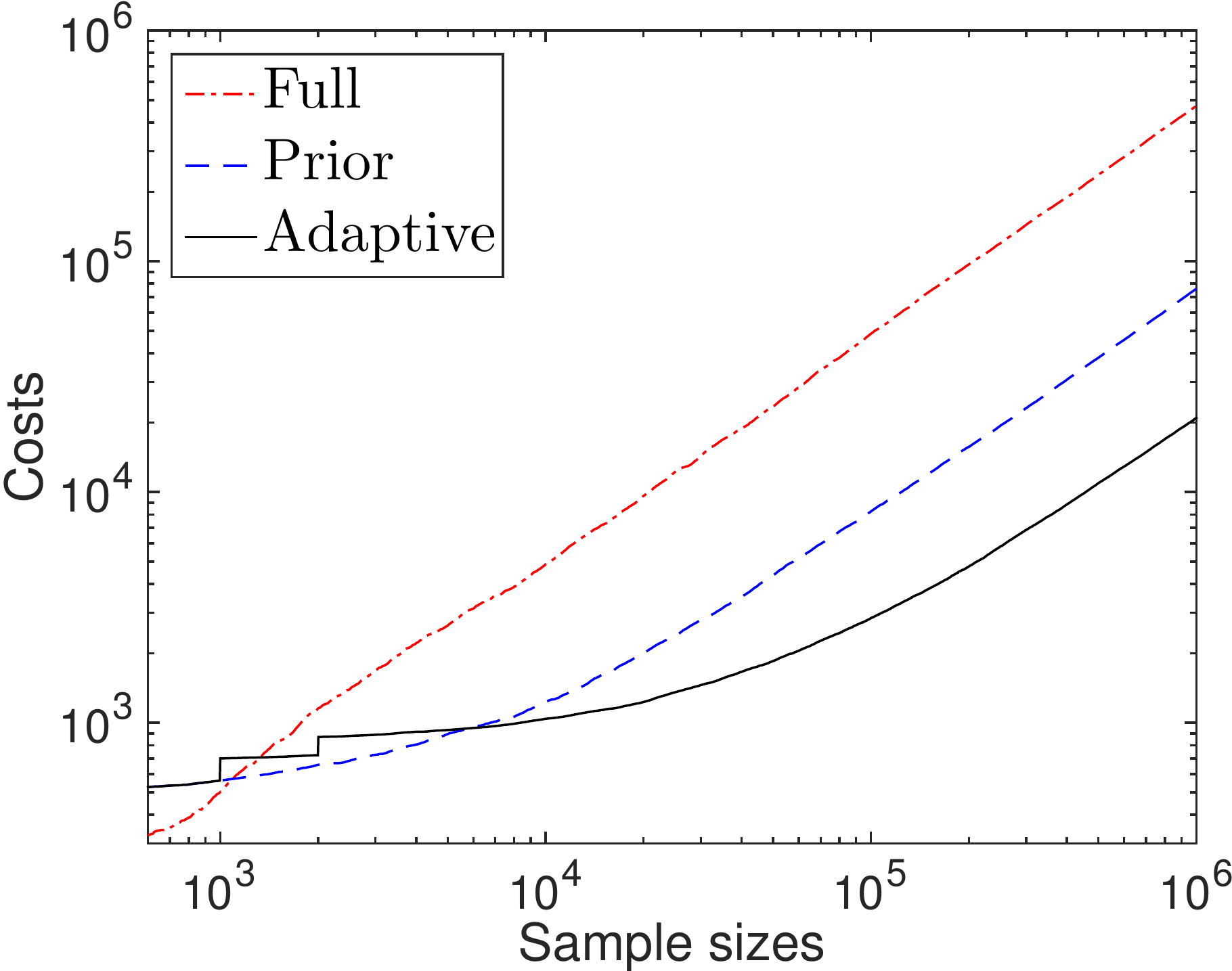} \\
(a) Costs for $M=23$ & (b) Costs for $M=73$ 
\end{tabular}}
\caption{Computational costs of full MCMC, prior RB-ANOVA-MCMC and adaptive RB-ANOVA-MCMC.}
\label{f:rb_cost}
\end{figure}

Figure~\ref{f:rb_cost} shows the costs with respect to the number of samples generated by the three methods. 
It is clear that, our adaptive RB-ANOVA-MCMC is the cheapest in the these three methods. 
From Figure~\ref{f:rb_cost}(a), to generate $10^6$ posterior samples for the test problem with $M=23$, 
the cost of adaptive RB-ANOVA-MCMC is around only one percent of the cost of full MCMC, and it is also much 
smaller than that of prior RB-ANOVA-MCMC. Note that the cost of full MCMC is slightly smaller than 
the sample size, since the prior distribution of the parameter $\xi$ is set to a uniform distribution in $[-1,1]^M$ and
the proposed samples are rejected without evaluating the forward model if they are outside of $[-1,1]^M$. 
 For the case of $M=73$ shown in Figure~\ref{f:rb_cost}(b), the cost 
of adaptive RB-ANOVA-MCMC to generate $10^6$ samples is around ten percent of full MCMC, and it is less than half of
the cost of prior RB-ANOVA-MCMC. From both Figure~\ref{f:rb_cost}(a) and Figure~\ref{f:rb_cost}(b), at an early stage
when the MCMC sample sizes are around $10^3$, adaptive RB-ANOVA-MCMC is more expensive than 
prior RB-ANOVA-MCMC. Moreover, for the case of $M=73$ shown in  Figure~\ref{f:rb_cost}(b), adaptive RB-ANOVA-MCMC is even 
more expensive than full MCMC at this early stage. The extra cost of adaptive RB-ANOVA-MCMC here comes from the reconstruction 
procedure (line 17 of Algorithm~\ref{alg_rbanovamcmc}). However, as the MCMC iteration continues, the reconstruction procedure 
quickly  stops, and the overall cost of adaptive RB-ANOVA-MCMC becomes much smaller than the costs of prior RB-ANOVA-MCMC and
full MCMC.

\begin{figure}[!htp]
\centerline{
\begin{tabular}{cc}
\includegraphics[width=5.5cm,height=4.5cm]{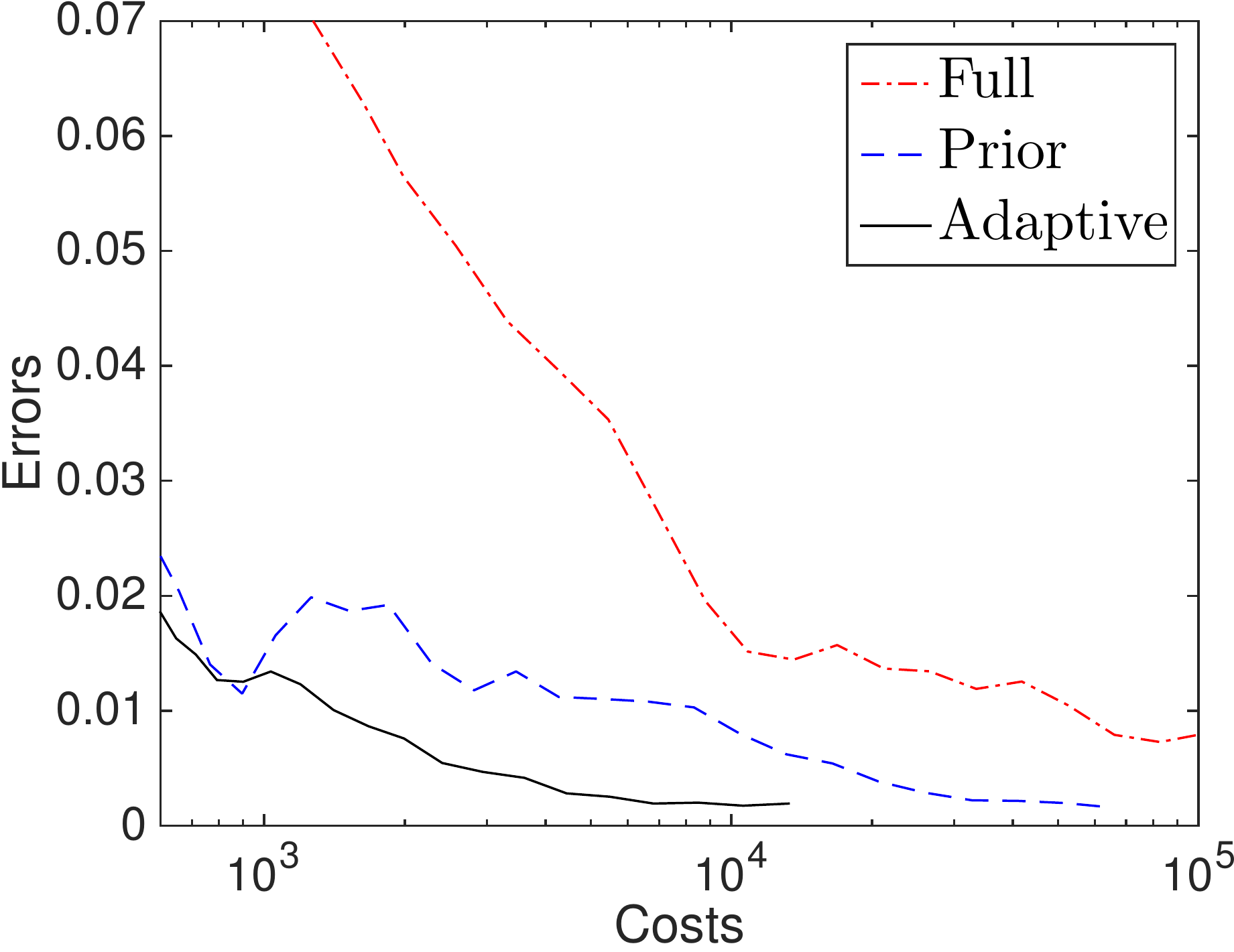} & 
\includegraphics[width=5.5cm,height=4.5cm]{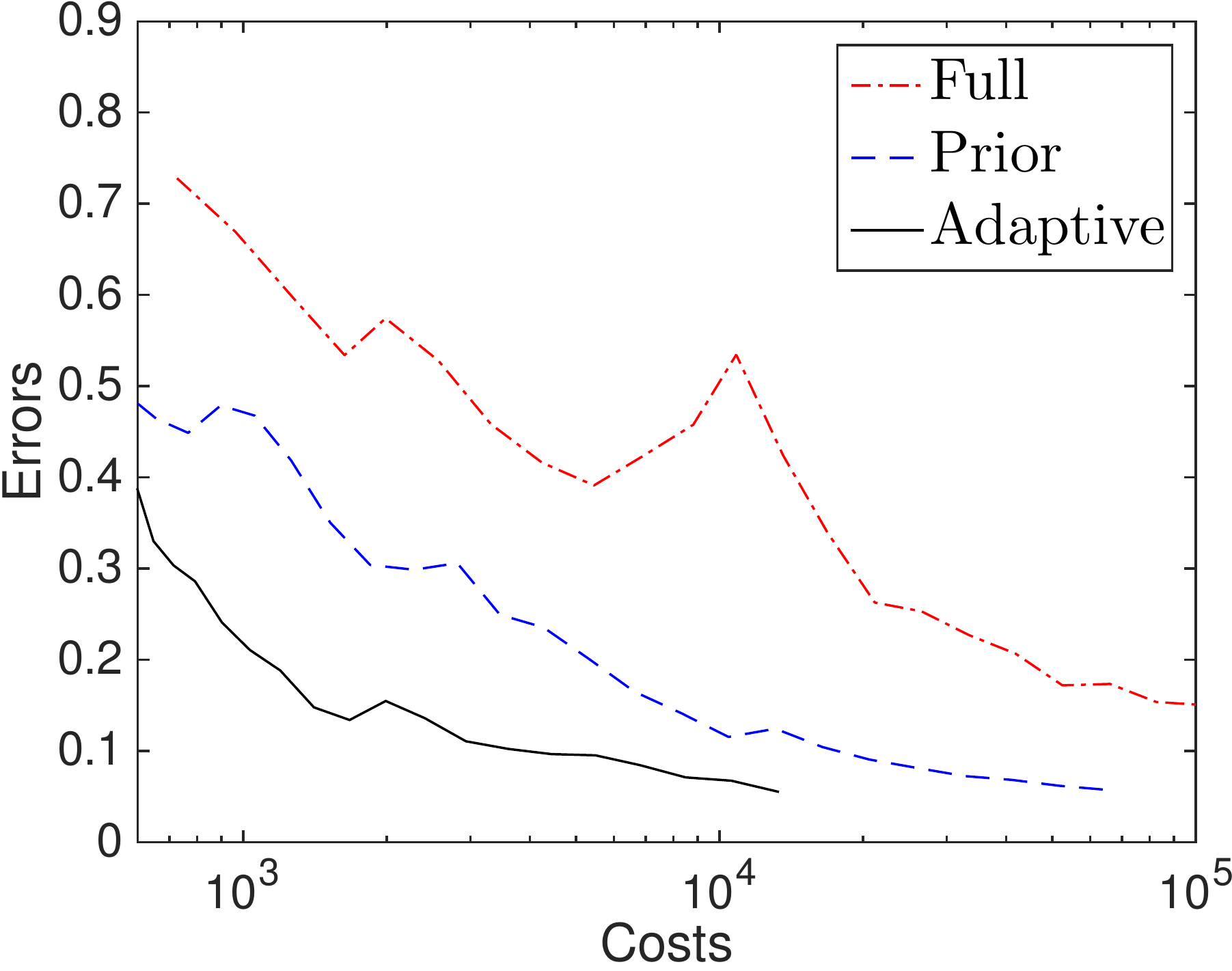} \\
(a) Mean errors & (b) Variance errors 
\end{tabular}}
\caption{Errors in mean and variance estimates ($\epsilon_{\rm mean}$ and $\epsilon_{\rm var}$) of full MCMC, prior RB-ANOVA-MCMC and adaptive RB-ANOVA-MCMC, 
for $\covl=5/4$ with $M=23$.}
\label{f:rb_errors_C125}
\end{figure}

\begin{figure}[!htp]
\centerline{
\begin{tabular}{cc}
\includegraphics[width=5.5cm,height=4.5cm]{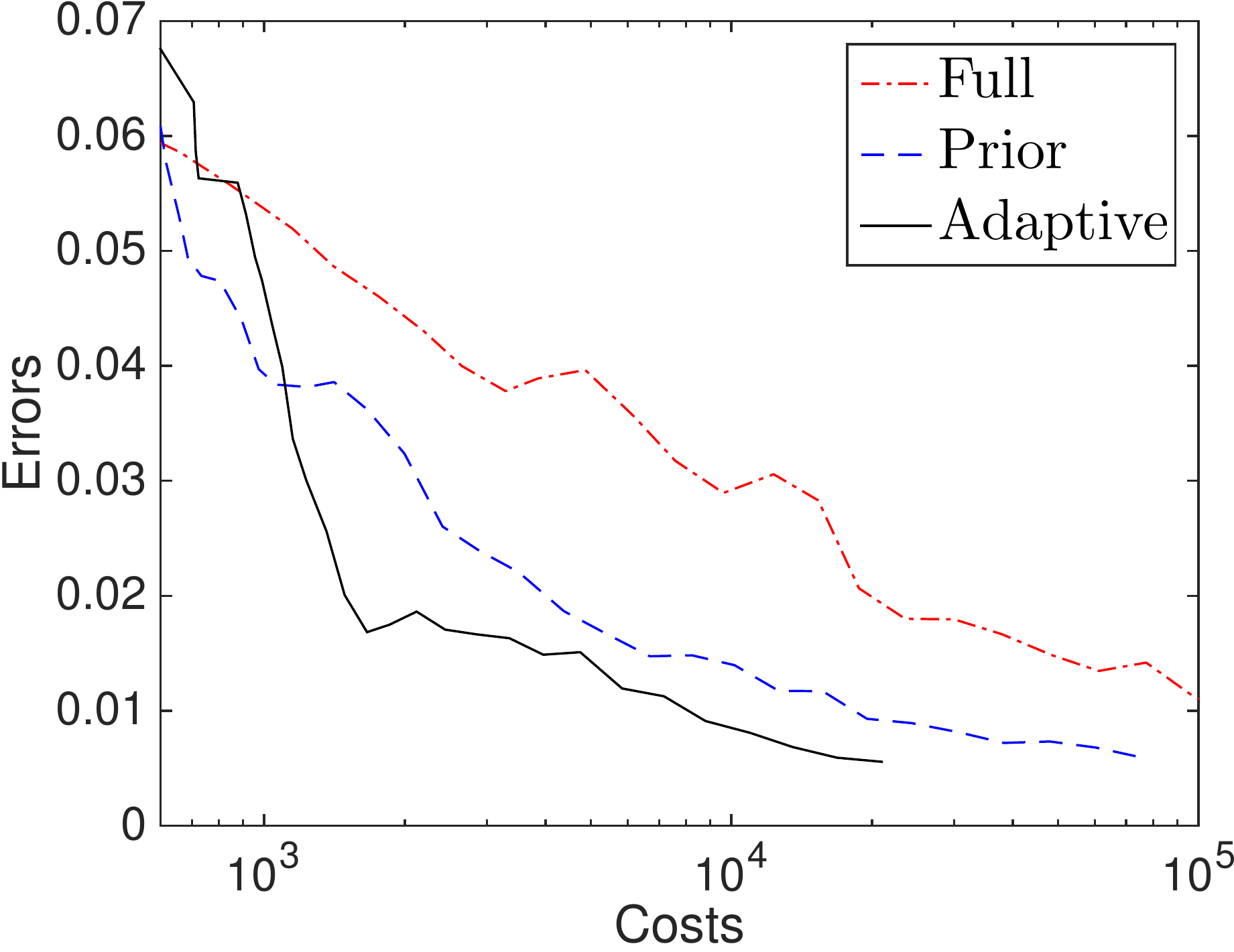} & 
\includegraphics[width=5.5cm,height=4.5cm]{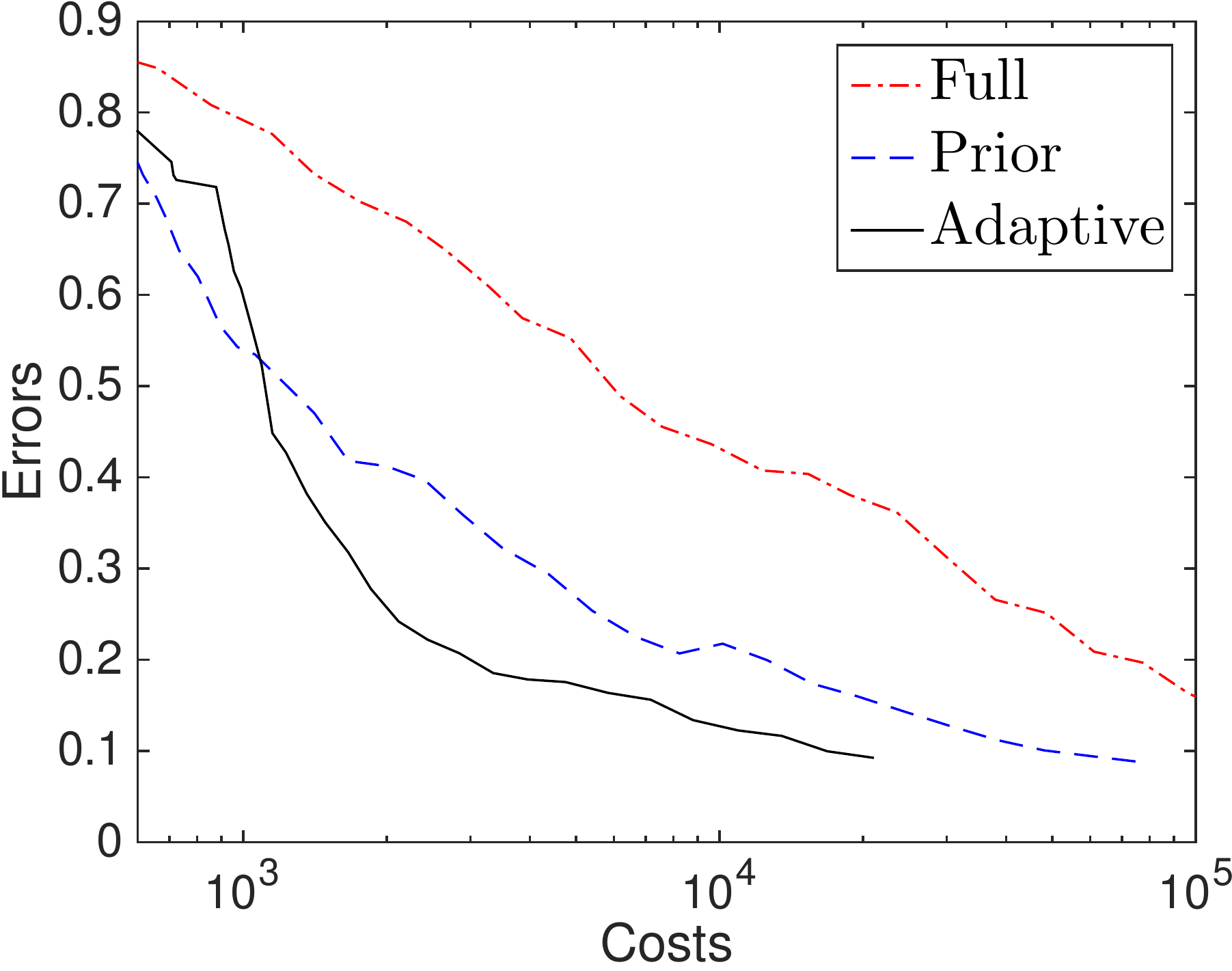} \\
(a) Mean errors & (b) Variance errors 
\end{tabular}}
\caption{Errors in mean and variance estimates ($\epsilon_{\rm mean}$ and $\epsilon_{\rm var}$) of full MCMC, prior RB-ANOVA-MCMC and adaptive RB-ANOVA-MCMC, 
for $\covl=5/8$ with $M=73$.}
\label{f:rb_errors_C0625}
\end{figure}

To assess the accuracy of RB-ANOVA-MCMC, 
we evaluate the errors in mean and variance estimates through the following quantities,  
\begin{subequations} \label{e:error}
\begin{align}
\epsilon_{\rm mean}:=\left\|\pE_{\Xi^*} \left(a\left(x,\xi\right) \right)-\pE_{\rm ref}\right\|_0 
\left/ 
\left\|\pE_{\rm ref}\right\|_0 \right., 
\label{e_m_anova}\\
\epsilon_{\rm var}:=\left\|\pV_{\Xi^*} \left(a\left(x,\xi\right) \right)-\pV_{\rm ref}\right\|_0 
\left/ 
\left\|\pV_{\rm ref}\right\|_0 \right.,
\label{e_m_ranova}
\end{align}
\end{subequations}
where $\pE_{\Xi^*}(a(x,\xi) )$ and $\pV_{\Xi^*}(a(x,\xi) )$ are defined in \eqref{mc_mean_field} and \eqref{mc_var_field}, 
and the reference mean estimate $\pE_{\rm ref}$ and the reference variance estimate $\pV_{\rm ref}$
are generated by full MCMC with $10^6$ samples using \eqref{mc_mean_field} and \eqref{mc_var_field}.  
Figure~\ref{f:rb_errors_C125} and Figure~\ref{f:rb_errors_C0625} show the errors of 
full MCMC, prior and adaptive RB-ANOVA-MCMC with respect to computational costs for the test problems 
with $M=23$ and $M=73$ respectively. It is clear that, the adaptive RB-ANOVA-MCMC method
has the smallest errors when the costs are not very small. 
For very small cost values, e.g., around $10^3$ in Figure~\ref{f:rb_errors_C0625}(a), 
the inefficiency of adaptive RB-ANOVA-MCMC (large errors in mean estimates) here is caused by the reconstruction procedure. As
the MCMC iteration continues and the cost values increase, cost spent in the reconstruction procedure of the adaptive approach 
becomes invisible, and the adaptively constructed model becomes significantly efficient . 
For example, for the case with $M=23$ shown in Figure~\ref{f:rb_errors_C125}(a), 
to achieve an accuracy in estimating the mean with error smaller than $0.01$,
the cost required by adaptive RB-ANOVA-MCMC is less than $2000$, 
which is less than a quarter of the cost required by prior RB-ANOVA-MCMC and is only around five percent of the cost of
required full MCMC. From Figure~\ref{f:rb_errors_C125}(b),
to achieve an accuracy in estimating the variance with error smaller than $0.2$ in this case, the cost of adaptive RB-ANOVA-MCMC 
is only around $1000$, which is only around twenty percent of the cost required by prior RB-ANOVA-MCMC
and is around two percent of the cost required by full RB-ANOVA-MCMC.
Similarly, for the case with $M=73$, Figure~\ref{f:rb_errors_C0625}(a) and Figure~\ref{f:rb_errors_C0625}(b)
show that to achieve given accuracies in mean and variance estimates, adaptive RB-ANOVA-MCMC  requires much less 
costs than  prior RB-ANOVA-MCMC and full MCMC.

Finally, the acceptance rates for generating $10^6$ posterior samples using the three approaches are shown in Table~\ref{table_acceptance_rate}.
It is clear that for both cases ($M=23$ and $M=73$), the acceptance rates of prior and adaptive RB-ANOVA-MCMC are consistent with 
the rates of full MCMC---they are the same up to two decimal places.

\begin{table}[!htp]
\caption{Acceptance rates of full MCMC, prior RB-ANOVA-MCMC and  adaptive RB-ANOVA-MCMC to generate $10^6$ posterior samples.}
\begin{center}
\renewcommand{\arraycolsep}{0.0pt}{
\begin{tabular}{l| c c c   c  }
\hline
$M$      &   Full    &   Prior   &  Adaptive  \\
$23$  &  $0.4175$ & $0.4193$ & $0.4159$ \\
$73$  & $0.2605$  &$0.2656$ & $0.2642$ \\
\hline
 \end{tabular}
}
\end{center}
\label{table_acceptance_rate}
\end{table}

\section{Conclusions}
\label{section_conclude}
Conducting posterior-oriented model reduction is one of the fundamental concepts for solving high-dimensional 
Bayesian inverse problems. With a focus on ANOVA, this paper proposes a novel adaptive reduced basis ANOVA 
(RB-ANOVA) model with respect to posterior distributions to accelerate MCMC procedures. 
The first novelty of our new approach is the adaptive ANOVA decomposition based on the posterior mean estimates.
It is known that the efficiency of the ANOVA decomposition is dependent on the choices of anchor points. 
Through adaptively updating the anchor point by posterior mean estimates during MCMC iterations,
an efficient ANOVA decomposition is obtained. 
Second, for all ANOVA terms, physical reduced bases are generated based on the posterior samples, 
which restricts the greedy algorithm to these samples so as to obtain optimal physical approximation bases for the Bayesian inversion 
problem. Numerical results demonstrate the overall efficiency of the proposed RB-ANOVA-MCMC algorithm.
As our algorithm is based on ANOVA decomposition with a single anchor point, it  currently can only be applied to
Bayesian inversion problems with unimodal posterior distributions. For multimodal distributions, a possible solution is  
to do ANOVA decomposition with multiple anchor points. Designing and analyzing ANOVA decomposition 
with multiple anchor points for both forward and inverse UQ problems will be the focus of our future work.

\bigskip
\textbf{Acknowledgments:}
Q. Liao  is support by NSFC under grant number 11601329
and J. Li is supported by the NSFC under grant number 11771289.

\section*{Reference}
\addcontentsline{toc}{chapter}{Bibliography}
\bibliography{liao,rb}

\begin{thebibliography}{10}
\expandafter\ifx\csname url\endcsname\relax
  \def\url#1{\texttt{#1}}\fi
\expandafter\ifx\csname urlprefix\endcsname\relax\def\urlprefix{URL }\fi
\expandafter\ifx\csname href\endcsname\relax
  \def\href#1#2{#2} \def\path#1{#1}\fi

\bibitem{tarantola2005inverse}
A.~Tarantola, Inverse problem theory and methods for model parameter
  estimation, SIAM, 2005.

\bibitem{tarantola2006popper}
A.~Tarantola, Popper, {B}ayes and the inverse problem, Nature physics 2~(8)
  (2006) 492--494.

\bibitem{kaipio2006statistical}
J.~Kaipio, E.~Somersalo, Statistical and computational inverse problems, Vol.
  160, Springer Science \& Business Media, 2006.

\bibitem{robert2005monte}
C.~P. Robert, G.~Casella, Monte carlo statistical methods (springer texts in
  statistics).

\bibitem{virieux2009overview}
J.~Virieux, S.~Operto, An overview of full-waveform inversion in exploration
  geophysics, Geophysics 74~(6) (2009) WCC1--WCC26.

\bibitem{yeh1986review}
W.~W.-G. Yeh, Review of parameter identification procedures in groundwater
  hydrology: The inverse problem, Water Resources Research 22~(2) (1986)
  95--108.

\bibitem{marzouk2009stochastic}
Y.~Marzouk, D.~Xiu, A stochastic collocation approach to {B}ayesian inference
  in inverse problems, Communications in Computational Physics 6~(4) (2009)
  826--847.

\bibitem{marzouk2009dimensionality}
Y.~M. Marzouk, H.~N. Najm, Dimensionality reduction and polynomial chaos
  acceleration of {B}ayesian inference in inverse problems, Journal of
  Computational Physics 228~(6) (2009) 1862--1902.

\bibitem{marzouk2007stochastic}
Y.~M. Marzouk, H.~N. Najm, L.~A. Rahn, Stochastic spectral methods for
  efficient {B}ayesian solution of inverse problems, Journal of Computational
  Physics 224~(2) (2007) 560--586.

\bibitem{nagel2016spectral}
J.~B. Nagel, B.~Sudret, Spectral likelihood expansions for {B}ayesian
  inference, Journal of Computational Physics 309 (2016) 267--294.

\bibitem{yan2015stochastic}
L.~Yan, L.~Guo, Stochastic collocation algorithms using l\_1-minimization for
  {B}ayesian solution of inverse problems, SIAM Journal on Scientific Computing
  37~(3) (2015) A1410--A1435.

\bibitem{bilionis2013solution}
I.~Bilionis, N.~Zabaras, Solution of inverse problems with limited forward
  solver evaluations: a {B}ayesian perspective, Inverse Problems 30~(1) (2013)
  015004.

\bibitem{kennedy2001bayesian}
M.~C. Kennedy, A.~O'Hagan, {B}ayesian calibration of computer models, Journal
  of the Royal Statistical Society: Series B (Statistical Methodology) 63~(3)
  (2001) 425--464.

\bibitem{wang2017adaptive}
H.~Wang, J.~Li, Adaptive {G}aussian process approximation for {B}ayesian
  inference with expensive likelihood functions, arXiv preprint
  arXiv:1703.09930.

\bibitem{ma2009efficient}
X.~Ma, N.~Zabaras, An efficient {B}ayesian inference approach to inverse
  problems based on an adaptive sparse grid collocation method, Inverse
  Problems 25~(3) (2009) 035013.

\bibitem{cui2015data}
T.~Cui, Y.~M. Marzouk, K.~E. Willcox, Data-driven model reduction for the
  {B}ayesian solution of inverse problems, International Journal for Numerical
  Methods in Engineering 102~(5) (2015) 966--990.

\bibitem{galbally2010non}
D.~Galbally, K.~Fidkowski, K.~Willcox, O.~Ghattas, Non-linear model reduction
  for uncertainty quantification in large-scale inverse problems, International
  journal for numerical methods in engineering 81~(12) (2010) 1581--1608.

\bibitem{lieberman2010parameter}
C.~Lieberman, K.~Willcox, O.~Ghattas, Parameter and state model reduction for
  large-scale statistical inverse problems, SIAM Journal on Scientific
  Computing 32~(5) (2010) 2523--2542.

\bibitem{wang2005using}
J.~Wang, N.~Zabaras, Using {B}ayesian statistics in the estimation of heat
  source in radiation, International Journal of Heat and Mass Transfer 48~(1)
  (2005) 15--29.

\bibitem{CMCS-CHAPTER-2008-001}
C.~Nguyen, G.~Rozza, D.~B.~P. Huynh, A.~T. Patera,
  \href{http://augustine.mit.edu/methodology/methodology_technical_papers.htm}{Reduced
  basis approximation and a posteriori error estimation for parametrized
  parabolic {PDE}s; {A}pplication to real-time {B}ayesian parameter
  estimation}, in: L.~Tenorio, B.~van Bloemen~Waanders, B.~Mallick, K.~Willcox,
  L.~Biegler, G.~Biros, O.~Ghattas, M.~Heinkenschloss, D.~Keyes, Y.~Marzouk
  (Eds.), Large {S}cale {I}nverse {P}roblems and {Q}uantification of
  {U}ncertainty, no. Chapter 8 in Wiley Series in Computational Statistics,
  John Wiley \& Sons, UK, 2010, pp. 151--178, ePFL-IACS report 11.2008.
\newline\urlprefix\url{http://augustine.mit.edu/methodology/methodology_technical_papers.htm}

\bibitem{frangos2010}
M.~Frangos, Y.~Marzouk, K.~Willcox, B.~van Bloemen~Waanders, Surrogate and
  Reduced-Order Modeling: A Comparison of Approaches for Large-Scale
  Statistical Inverse Problems, John Wiley \& Sons, Ltd, 2010, pp. 123--149.

\bibitem{li2015note}
J.~Li, A note on the {K}arhunen--{L}o{\`e}ve expansions for
  infinite-dimensional {B}ayesian inverse problems, Statistics \& Probability
  Letters 106 (2015) 1--4.

\bibitem{fisher}
R.~Fisher, Statistical Methods for Research Workers, Oliver and Boyd, Berlin,
  1925.

\bibitem{sob03}
I.~Sobol, Theorems and examples on high dimensional model representation,
  Reliability Engineering and System Safety 79 (2003) 187--193.

\bibitem{caogun09}
Y.~Cao, Z.~Chen, M.~Gunzburger, {ANOVA} expansions and efficient sampling
  methods for parameter dependent nonlinear {PDE}s, International Journal of
  Numerical Analysis and Modeling 6 (2009) 256--273.

\bibitem{wintar09}
C.~Winter, A.~Guadagnini, D.~Nychka, D.~Tartakovsky, Multivariate sensitivity
  analysis of saturated flow through simulated highly heterogeneous groundwater
  aquifers, Journal of Computational Physics 217 (2009) 166--175.

\bibitem{gaohes10}
Z.~Gao, J.~S. Hesthaven, On {ANOVA} expansions and strategies for choosing the
  anchor point, Applied Mathematics and Computation 217 (2010) 3274--3285.

\bibitem{zhachokar11}
Z.~Zhang, M.~Choi, G.~Karniadakis, Anchor points matter in {ANOVA}
  decomposition, Spectral and High Order Methods for Partial Diferential
  Equations Lecture Notes in Computational Science and Engineering 76 (2011)
  347--355.

\bibitem{maza10}
X.~Ma, N.~Zabaras, An adaptive high-dimensional stochastic model representation
  technique for the solution of stochastic partial differential equations,
  Journal of Computational Physics 229 (2010) 3884--3915.

\bibitem{yanlin11}
X.~Yang, M.~Choi, G.~Lin, G.~E. Karniadakis, Adaptive {ANOVA} decomposition of
  stochastic incompressible and compressible flows, Journal of Computational
  Physics 231 (2012) 1587--1614.

\bibitem{fookar09}
J.~Foo, G.~Karniadakis, Multi-element probabilistic collocation in high
  dimensions, Journal of Computational Physics 229 (2010) 1536--1557.

\bibitem{zhahes16}
J.~S. Hesthaven, S.~Zhang, On the use of {ANOVA} expansions in reduced basis
  methods for high-dimensional parametric partial differential equations,
  Journal of Scientific Computing, {}To appear, DOI: 10.1007/s10915-016-0194-9.

\bibitem{liaolin16}
Q.~Liao, G.~Lin, Reduced basis {ANOVA} methods for partial differential
  equations with high-dimensional random inputs, Journal of Computational
  Physics 317 (2016) 148--164.

\bibitem{choelman18}
H.~Cho, H.~C. Elman, An adaptive reduced basis collocation method based on
  {PCM} {ANOVA} decomposition for anisotropic stochastic {PDEs}, International
  Journal for Uncertainty Quantification 8 (2018) 193--210.

\bibitem{zhachokar12}
Z.~Zhang, M.~Choi, G.~Karniadakis, Error estimates for the {ANOVA} method with
  polynomial chaos interpolation: Tensor product functions, SIAM Journal on
  Scientific Computing 34~(2) (2012) A1165--A1186.

\bibitem{sirovich87}
L.~Sirovich, Turbulence and the dynamics of coherent structures, {P}art {I}:
  Coherent structures, Quarterly of Applied Mathematics 45 (1987) 561--571.

\bibitem{holmes96}
P.~Holmes, J.~L. Lumley, G.~Berkooz, Turbulence, Coherent Structures, Dynamical
  Systems and Symmetry, Cambridge, New York, 1996.

\bibitem{gunpetsha07}
M.~Gunzburger, J.~Peterson, J.~Shadid, Reduced-order modeling of time-dependent
  {PDE}s with multiple parameters in the boundary data, Computer Methods in
  Applied Mechanics and Engineering 196 (2007) 1030--1047.

\bibitem{verpat02}
K.~Veroy, D.~Rovas, A.~Patera, A posteriori error estimation for reduced-basis
  approximation of parametrized elliptic coercive partial differential
  equations:{``Convex Inverse"} bound conditioners, ESAIM: Control,
  Optimisation and Calculus of Variations 8 (2002) 1007--1028.

\bibitem{nguver05}
N.~Nguyen, K.~Veroy, A.~Patera, Certified real-time solution of parametrized
  partial differential equations, in: S.~Yip (Ed.), Handbook of Materials
  Modeling, Springer, 2005, pp. 1523--1558.

\bibitem{haaohl08}
B.~Haasdonk, M.~Ohlberger, Reduced basis method for finite volume
  approximations of parametrized linear evolution equations, ESAIM:
  Mathematical Modelling and Numerical Analysis 42 (2008) 277--302.

\bibitem{buiwillcox08}
T.~Bui-Thanh, K.~Willcox, O.~Ghattas, Model reduction for large-scale systems
  with high-dimensional parametric input space, SIAM Journal on Scientific
  Computing 30 (2008) 3270--3288.

\bibitem{boybri10}
S.~Boyaval, C.~L. Bris, T.~Leli\`evre, Y.~Maday, N.~Nguyen, A.~Patera, Reduced
  basis techniques for stochastic problems, Archives of Computational Methods
  in Engineering 17 (2010) 1--20.

\bibitem{patrozbook}
A.~Patera, G.~Rozza, Reduced Basis Approximation and A Posteriori Error
  Estimation for Parametrized Partial Differential Equations, 2007, version
  1.0, Copyright MIT 2006--2007, to appear in (tentative title) MIT Pappalardo
  Graduate Monographs in Mechanical Engineering.

\bibitem{quamanbook}
A.~Quarteroni, A.~Manzoni, F.~Negri, Reduced Basis Methods for Partial
  Differential Equations, Springer International Publishing, Springer
  International Publishing Switzerland, 2016.

\bibitem{elmanliao}
H.~Elman, Q.~Liao, Reduced basis collocation methods for partial differential
  equations with random coefficients, SIAM/ASA Journal on Uncertainty
  Quantification 1 (2013) 192--217.

\bibitem{powellnewsum16}
C.~Newsum, C.~Powell, Efficient reduced basis methods for saddle point problems
  with applications in groundwater flow, SIAM/ASA Journal on Uncertainty
  Quantification 5~(1) (2017) 1248--1278.

\bibitem{zhangwebster17}
Q.~Guan, M.~Gunzburger, C.~G. Webster, G.~Zhang, Reduced basis methods for
  nonlocal diffusion problems with random input data, Computer Methods in
  Applied Mechanics and Engineering 317 (2017) 746 -- 770.

\bibitem{li2014adaptive}
J.~Li, Y.~M. Marzouk, Adaptive construction of surrogates for the {Bayesian}
  solution of inverse problems, SIAM Journal on Scientific Computing 36~(3)
  (2014) A1163--A1186.

\bibitem{brae97}
D.~Braess, Finite Elements, Cambridge University Press, London, 1997.

\bibitem{elman05}
H.~Elman, D.~Silvester, A.~Wathen, Finite Elements and Fast Iterative Solvers,
  Oxford University Press, New York, 2005.

\bibitem{ghaspa03}
R.~Ghanem, P.~Spanos, Stochastic Finite Elements: A Spectral Approach, Dover
  Publications, New York, 2003.

\bibitem{babuska1}
I.~Babu\v{s}ka, F.~Nobile, R.~Tempone, A stochastic collocation method for
  elliptic partial differential equations with random input data, SIAM Journal
  on Numerical Analysis 45 (2007) 1005--1034.

\bibitem{elmmil11}
H.~Elman, C.~Miller, E.~Phipps, R.~Tuminaro, Assessment of collocation and
  {G}alerkin approaches to linear diffusion equations with random data,
  International Journal for Uncertainty Quantification 1 (2011) 19--34.

\bibitem{powelm09}
C.~Powell, H.~Elman, Block-diagonal preconditioning for spectral stochastic
  finite-element systems, IMA Journal of Numerical Analysis 29 (2009) 350--375.

\bibitem{roror01}
G.~O. Roberts, J.~S. Rosenthal, Optimal scaling for various metropolis-hastings
  algorithms, Statistical Science 16 (2001) 351--367.

\end{thebibliography}

\end{document}